\documentclass[12pt]{article} %fleqn
 \usepackage{ams}
 \usepackage{dsfont}
 \usepackage{amsmath}
 \usepackage{amsthm}
 \usepackage{txfonts,bm}
 \usepackage{bold-extra}
 \usepackage[sc,bf]{caption}
 \usepackage{rotating}
 \usepackage{fancyhdr}
 \usepackage{fancyheadings}
 \RequirePackage[OT1]{fontenc}
 \RequirePackage{amsthm,amsmath}
 \RequirePackage[numbers]{natbib}
 \usepackage[colorlinks,citecolor=blue,urlcolor=blue,linkcolor=blue]{hyperref}

 \usepackage{titlesec}
 \usepackage{latexsym}
   \font\twelvebm                       = cmmib10 at 12truept
   \font\tenbm                          = cmmib10 at 10truept
   \font\sevenbm                        = cmmib10 at 7truept
 = \tenbm \scriptscriptfont9              = \sevenbm
 
   at 10truept  at 10truept  at
 10truept  at 10truept
  at 10truept
 \mathchardef \BGamma            = "0900 \mathchardef \BDelta
 = "0901 \mathchardef \BTheta            = "0902 \mathchardef
 \BLambda           = "0903 \mathchardef \BXi               = "0904
 \mathchardef \BPi               = "0905 \mathchardef \BSigma
 = "0906 \mathchardef \BUpsilon          = "0907 \mathchardef \BPhi
 = "0908 \mathchardef \BPsi              = "0909 \mathchardef
 \BOmega            = "090A \mathchardef \Balpha            = "090B
 \mathchardef \Bbeta             = "090C \mathchardef \Bgamma
 = "090D \mathchardef \Bdelta            = "090E \mathchardef
 \Bepsilon          = "090F \mathchardef \Bzeta             = "0910
 \mathchardef \Beta              = "0911 \mathchardef \Btheta
 = "0912 \mathchardef \Biota             = "0913 \mathchardef
 \Bkappa            = "0914 \mathchardef \Blambda           = "0915
 \mathchardef \Bmu               = "0916 \mathchardef \Bnu
 = "0917 \mathchardef \Bxi               = "0918 \mathchardef \Bpi
 = "0919 \mathchardef \Brho              = "091A \mathchardef
 \Bsigma            = "091B \mathchardef \Btau              = "091C
 \mathchardef \Bupsilon          = "091D \mathchardef \Bphi
 = "091E \mathchardef \Bchi              = "091F \mathchardef \Bpsi
 = "0920 \mathchardef \Bomega            = "0921 \mathchardef
 \Bvarepsilon       = "0922 \mathchardef \Bvartheta         = "0923
 \mathchardef \Bvarpi            = "0924 \mathchardef \Bvarrho
 = "0925 \mathchardef \Bvarsigma         = "0926 \mathchardef
 \Bvarphi           = "0927
 \mathchardef \bA        = "0941 \mathchardef \bB        = "0942
 \mathchardef \bC        = "0943 \mathchardef \bD        = "0944
 \mathchardef \bE        = "0945 \mathchardef \bF        = "0946
 \mathchardef \bG        = "0947 \mathchardef \bH        = "0948
 \mathchardef \bI        = "0949 \mathchardef \bJ        = "094A
 \mathchardef \bK        = "094B \mathchardef \bL        = "094C
 \mathchardef \bM        = "094D \mathchardef \bN        = "094E
 \mathchardef \bO        = "094F \mathchardef \bP        = "0950
 \mathchardef \bQ        = "0951 \mathchardef \bR        = "0952
 \mathchardef \bS        = "0953 \mathchardef \bT        = "0954
 \mathchardef \bU        = "0955 \mathchardef \bV        = "0956
 \mathchardef \bW        = "0957 \mathchardef \bX        = "0958
 \mathchardef \bY        = "0959 \mathchardef \bZ        = "095A
 \mathchardef \ba        = "0961 \mathchardef \bb        = "0962
 \mathchardef \bc        = "0963 \mathchardef \bd        = "0964
 \mathchardef \bee       = "0965 %%%I CHANGED IT FROM \be; SEE IN THE BEGGINING.
 \mathchardef \bff       = "0966 \mathchardef \bg        = "0967
 \mathchardef \bh        = "0968
 \mathchardef \bj        = "096A \mathchardef \bk        = "096B
 \mathchardef \bl        = "096C \mathchardef \bm        = "096D
 \mathchardef \bn        = "096E \mathchardef \bo        = "096F
 \mathchardef \bp        = "0970 \mathchardef \bq        = "0971
 \mathchardef \br        = "0972 \mathchardef \bs        = "0973
 \mathchardef \bt        = "0974 \mathchardef \bu        = "0975
 \mathchardef \bv        = "0976 \mathchardef \bw        = "0977
 \mathchardef \bx        = "0978 \mathchardef \by        = "0979
 \mathchardef \bz        = "097A

 \font\tencb            = cmssbx10 scaled \magstep4 \font\eigcb
 = cmssbx10 scaled \magstep2 \textfont8             = \tencb
 \mathchardef\bAs       = "1841
 \def\Asem#1#2{\mathop{\vrule height10.5pt depth5.5pt width0pt\bAs}_{#1}^{#2}}
 \def\asem#1#2{
          \ifmmode
         \ifinner
            \raise0.9pt\hbox{$\scriptstyle\bAs$}_{#1}^{#2}
         \else
            \Asem{#1}{#2}
         \fi
          \fi
          }

 \newtheorem{theo}{\small\bf Theorem}[section]
 \newtheorem{lem}{\small\bf Lemma}[section]
 \newtheorem{prop}{\small\bf Proposition}[section]
 \newtheorem{rem}{\small\bf Remark}[section]
 
 \newtheorem{exam}{\small\bf Example}[section]
 
  \newtheorem{defi}{\small\bf Definition}[section]
 
 \newtheorem{cor}{\small\bf Corollary}[section]

 \newcommand{\be}{\begin{equation}}
 \newcommand{\ee}{\end{equation}}
 
 \newcommand{\lawv}{\stackrel{\mbox{\tiny d}}{\rightarrow}}

 \newcommand{\E}{\mbox{\rm \hspace*{.2ex}I\hspace{-.5ex}E\hspace*{.2ex}}}

 \newcommand{\Var}{\mbox{\rm \hspace*{.2ex}Var\hspace*{.2ex}}}
 
 \newcommand{\Cov}{\mbox{\rm \hspace*{.2ex}Cov\hspace*{.2ex}}}

 \newcommand{\bbb}[1]{\mbox{\boldmath $ #1 $}}

 \newenvironment{pr}[1]{{\small\bf {#1}:}}{}

 \newcommand{\ud}{\hspace{0.1ex}\textrm{$\rm{d}$}}
 \newcommand{\RR}{\mathbb{R}}
 
 \newcommand{\QQ}{\mathbb{Q}}

 \numberwithin{equation}{section}

 \title{\Large\bf On point estimators for
 Gamma and Beta distributions}
 \author{\large
 Nickos
 Papadatos%\textcolor{blue}{$\rm^b$}$^,$\textcolor[rgb]{0,0,1}
 {\footnote{
 {e-mail:}\
 \textcolor[rgb]{0.98,0.00,0.00}{npapadat@math.uoa.gr}, {url:}\
 \textcolor[rgb]{0.98,0.00,0.00}{users.uoa.gr/$\sim$npapadat/}}}}
 
 \vspace{-1em}

 \date{\small\it
 \begin{tabular}{r@{\hspace{0ex}}l}
 & National and Kapodistrian University of Athens, Department of Mathematics,
 \\
 [-.3ex]
 &
 Section of Statistics and Operations
 Research,
 Panepistemiopolis, 157 84 Athens, Greece.
 \end{tabular}
 }
 \vspace{-1em}

\begin{document}

 \maketitle
 \vspace*{-2em}

 \vspace{-12em}
 \noindent
 {\sc\hfill Dedicated to Professor Stavros Kourouklis}
 \vspace{12em}

 \thispagestyle{empty}

 \begin{abstract}
 %\vspace*{-.5em}
 \noindent
 Let $X_1,\ldots,X_n$ be a random sample from the Gamma distribution
 with density $f(x)=\lambda^{\alpha}x^{\alpha-1}e^{-\lambda x}/\Gamma(\alpha)$,
 $x>0$, where both $\alpha>0$ (the shape parameter) and $\lambda>0$
 (the reciprocal scale parameter) are unknown.
 The main result shows that the {\it uniformly minimum variance unbiased
 estimator} (UMVUE) of the shape parameter, $\alpha$, exists
 if and only if $n\geq 4$; moreover, it
 has finite variance if and only if $n\geq 6$.
 More precisely, the form of the UMVUE is given for all parametric
 functions $\alpha$, $\lambda$, $1/\alpha$ and $1/\lambda$.
 Furthermore, a highly efficient estimating  procedure for the two-parameter
 Beta distribution is also given. This is based on a Stein-type covariance
 identity for the Beta distribution, followed by an application
 of the theory of $U$-statistics and the delta-method.
 \end{abstract}
 {\footnotesize {\it MSC}:  Primary 62F10; 62F12; Secondary 62E15.
 \vspace*{-.7ex}
 Ye-Chen-type closed-form estimators;
 asymptotic efficiency; $U$-statistics; Stein-type covariance identity;
 delta-method.
 %; moment problem.
 }
 %\vspace*{.2em}
 %\newline
 %{\it RUNNING HEAD}: Factorial Moment Distance.}
 \vspace*{-1em}

 \section{Introduction}
 %\vspace*{-.5em}
 \setcounter{equation}{0}
 \label{sec.1}

 Let $\Gamma(\alpha)=\int_0^\infty x^{\alpha-1}e^{-x} \ud x$, $\alpha>0$,
 be the Gamma function.
 The Gamma distribution, $\mathcal{G}(\alpha,\lambda)$, with density
 \[
 f(x)=\frac{\lambda^\alpha}{\Gamma(\alpha)} x^{\alpha-1} e^{-\lambda x},
  \ \ \ \ x>0,
 \]
 and parameters $\alpha>0$ (the shape) and $\lambda>0$ (the reciprocal scale),
 generates one of the most useful statistical models, especially when the data
 are nonnegative. This distribution belongs to the Pearson family, and in particular,
 to the Integrated Pearson family of distributions, see, e.g.,
 Afendras and Papadatos (2015).
 Therefore, estimating procedures for the parameters $\alpha$, $\lambda$, and/or
 $1/\alpha$ and $1/\lambda$, are of fundamental importance in applications.
 It is well-known that the maximum likelihood estimators (MLEs) are not tractable and, consequently, one has to apply numerical procedures for their evaluation. Moreover,
 to the best of our knowledge, it is not even known whether
 an unbiased estimator for $\alpha$ (and $\lambda$) exists, for some
 sample size $n$. In the contrary, unbiased estimators
 for the reciprocal parameters $1/\alpha$ and $1/\lambda$ do exist;
 recently, Ye and Chen (2017) obtained closed-form unbiased estimators
 for the reciprocal parameters with very high asymptotic efficiency.
 It should be noted that an augmented-sample (simulation)
 technique for obtaining exact confidence intervals
 for $\alpha$ is introduced by Iliopoulos (2016).

 In the present work we show that an unbiased estimator for $\alpha$
 exists if and only if the sample size $n$ is at least $4$, and the same
 is true for $\lambda$ (provided that $\alpha$ is not too small).
 More precisely, in Section \ref{sec.2}
 we obtain the uniformly minimum variance unbiased estimators
 (UMVUEs) of $\alpha$, $\lambda$, $1/\alpha$, $1/\lambda$. In particular,
 the UMVUE of $1/\lambda$, based on two observations, can be written as a positive
 symmetric kernel, whence, the corresponding $U$-statistic
 generates the Ye-Chen (2017) estimator; see Remark \ref{rem.closed.form}.

 Let $B(\alpha,\beta)=\Gamma(\alpha)\Gamma(\beta)/\Gamma(\alpha+\beta)$, $\alpha>0$,
 $\beta>0$,  be the Beta function. The Beta distribution has density
 \be
 \label{beta.distribution}
 f(x)=\frac{1}{B(\alpha,\beta)} x^{\alpha-1} (1-x)^{\beta-1},
  \ \ \ \ 0<x<1.
 \ee
 This distribution belongs to the Integrated Pearson family,
 and, at least to our knowledge,
 it is not known whether the parameters can be unbiasedly estimated.
 The MLEs do not have closed forms, and one has to overcome similar difficulties
 as for the Gamma case.
 In Section \ref{sec.3} we obtain closed-form, Ye-Chen-type, estimators
 for the parameters.
 Our method is different than that of Ye and Chen (2017), since it is based on
 a Stein-type covariance identity for the Beta distribution, followed
 by an application of the theory of $U$-statistics. The covariance identity
 generates an unbiased bivariate symmetric kernel for the parameter
 $\theta=1/(\alpha+\beta)$.
 Then, an application of the classical theory of $U$-statistics
 produces an unbiased estimator of $\theta$,
 which, in turn, yields the desired estimators of $\alpha$ and $\beta$.
 The asymptotic
 efficiency of the proposed estimators
 is studied in some detail, and it is shown to be very high
 (see Table 1). The Beta model is useful not only for analyzing data
 with bounded support, but also when positive data are
 mapped to the interval $(0,1)$ by means of a transformation,
 e.g., $x\mapsto x/(1+x)$.

 \section{Unbiased estimation of the parameters of Gamma distribution}
 \setcounter{equation}{0}
 \label{sec.2}
 Let ${\bbb X}=(X_1,\ldots,X_n)$ ($n\geq 2$)
 be a random sample from ${\it \Gamma}(\alpha,\lambda)$,
 where the unknown parameter $\bbb{\theta}=(\alpha,\lambda)\in
 {\bbb \Theta}=(0,\infty)^2$.
 It is well-known that the complete, sufficient statistic
 can be written in the form
 ${\bbb T}(\bbb X)=
 %\Big(\overline{X}_n,(\overline{X}_n)^{-n}\prod_{i=1}^n X_i\Big)=
 %\Big(\overline{X}_n,Y_n\Big)$, where $\overline{X}_n=n^{-1}\sum_{i=1}^n X_i$
 (X,Y)$,  where $X=\sum_{i=1}^n X_i$
 and $Y^{1/n}$ is the ratio of the geometric to the arithmetic
 mean of the data. The random variables %$\overline{X}_n$ and $Y_n$
 $X,Y$ are independent, as follows, e.g.,
 from Basu's theorem, because the distribution of $Y$ is free of $\lambda$,
 while $X$ is complete for $\lambda$ when $\alpha>0$ is known.
 The random  variable $X$ follows a
 $\mathcal{G}(n\alpha,\lambda)$ distribution, while $Y$
 is distributed like $\prod_{i=1}^{n-1}B_\alpha(i/n)$,
 where the random variables $B_\alpha(i/n)$, $i=1,\ldots,n-1$, are
 independent, and $B_\alpha(i/n)$
 has a Beta distribution (\ref{beta.distribution}) with parameters $\alpha$ and
 $\beta=i/n$.
 %,
 %$B(\alpha,i/n)$, with density
 %$\beta(x;\alpha,i/n)
 %=\Gamma(\alpha+i/n) x^{\alpha-1}(1-x)^{i/n-1}/(\Gamma(\alpha)\Gamma(i/n))$,
 %$0<x<1$.

 By using the Beta product representation,
 it follows inductively that the random variable $Y$ has density
 \be
 \label{f_Y}
 f_Y(y;\alpha)=K(\alpha) y^{\alpha-1}(1-y)^{(n-3)/2}
 G(y),
 \ \ \
 0<y<1,
 \ee
 where
 $K(\alpha)=\Gamma(\alpha)^{-n}\prod_{i=1}^{n-1}\Gamma(\alpha+i/n)/\Gamma(i/n)$,
 and the (power-series) function $G$
 is given by
 \be
 \label{g_(n-1)}
 G(y)=\sum_{m=0}^\infty d_{m}(1-y)^m, \ \ 0<y<1.
 \ee
 Notice that our notation suspends the dependence of $K(\alpha)$
 and $G(y)$ on $n$.

 \noindent
 In (\ref{g_(n-1)}) the coefficients $(d_{m})_{m\geq 0}$ are defined
 recurrently as follows:
 \begin{eqnarray}
 \nonumber
 &&
 \gamma_1(m)=I(m=0)=\left\{\begin{array}{ccl} 1, & \mbox{if} & m=0,\\
       0,& \mbox{if} & m=1,2,\ldots,  \end{array}\right.
       \\
       \nonumber
       &&
       \gamma_{i+1}(m)=\frac{\Gamma(m+i(i+1)/(2n))}{\Gamma(m+(i+1)(i+2)/(2n))}
       \sum_{k=0}^m \frac{\Gamma(m-k+(i+1)/n)}{(m-k)!}\gamma_i(k),
       \ \ \ \ i=1,\ldots,n-2,
       \\
       &&
       \label{gamma}
       d_m:=\gamma_{n-1}(m), \ \ \ \ m=0,1,\ldots \ . \Bigg.
 \end{eqnarray}
 Most of the preceding results can be found in Glaser (1976), Nandi (1980),
 and Tang and Gupta (1984). Obviously,
 $G(y)\equiv 1$ when $n=2$.
 Moreover, for $n\geq 3$,
 $G$ is strictly decreasing and positive; more precisely,
 $d_m>0$ for all $m$, $G'(y)<0$ and $G(y)>G(1-)=d_0>0$, $0<y<1$.
 It should be noted at this point that, when $n\geq 3$,
 the auxiliary random variable
 $Y_1=\prod_{i=1}^{n-2}B_\alpha(i/n)$ has density
 (cf.\ Tang and Gupta, 1984)
 \be
 \label{f_Y1}
 f_{Y_1}(y;\alpha)=K_1(\alpha) y^{\alpha-1}(1-y)^{(n-1)(n-2)/(2n)-1} G_1(y),
 \ \ \ 0<y<1,
 \ee
 where $K_1(\alpha)=\Gamma(\alpha)^{-(n-1)}\prod_{i=1}^{n-2}
 \Gamma(\alpha+i/n)/\Gamma(i/n)$ and
 \be
 \label{g1}
 G_1(y)=\sum_{k=0}^\infty \gamma_{n-2}(k)(1-y)^k, \ \ \ 0<y<1.
 \ee
 The random variable $Y_1$
 will be used in the proof of Theorem \ref{theo.unbiased}, below.

 \begin{defi}{\rm
 \label{def.UMVUE}
 An estimator $T(\bbb X)$ will be called {\it uniformly minimum variance
 unbiased estimator} (UMVUE) for the parametric function
 $h(\alpha,\lambda)$ if $T(\bbb X)=u(X,Y)$ is a (Borel)
 function of the complete, sufficient statistic $(X,Y)$, and
 $\E_{(\alpha,\lambda)}u(X,Y)=h(\alpha,\lambda)$, for all
 $\alpha>0$, $\lambda>0$, where the subscript $(\alpha,\lambda)$
 denotes expectation w.r.t.\ the joint density of
 $(X_1,\ldots,X_n)$, %namely,
 \[
 f(x_1,\ldots,x_n)=\frac{\lambda^{n\alpha}}{\Gamma(\alpha)^n}
 \left(\prod_{i=1}^n x_i\right)^{\alpha-1} \exp\left(-\lambda\sum_{i=1}^n x_i\right),
 \ \ \ \ (x_1,\ldots,x_n)\in(0,\infty)^n.
 \]
 Accordingly, an UMVUE may, or may not have
 finite variance. For completeness of the presentation,
 the term UMVUE will be also used
 in a wider sense, namely, in the case  where
 %the function $u(x,y)$
 %is the restriction of an holomorphic function
 %$u(z_1,z_2):\CC_+\times \CC_2\to \CC$,
 %where $\CC_+=\{z\in\CC:\Re(z)>0\}$, $\CC_2=\{z\in\CC:0<\Re(z)<1\}$,
 %and
 the equality
 $\E_{(\alpha,\lambda)}u(X,Y)=h(\alpha,\lambda)$
 is merely satisfied for all $(\alpha,\lambda)$ in
 a subset of $(0,\infty)^2$ with nonempty interior.
 %Here, the term
 %holomorphic means that for each fixed $z_2\in\CC_2$,
 %$u_1(z):=u(z,z_2)$ is differentiable (in the complex sense) for all $z\in \CC_+$,
 %and for each fixed $z_1\in\CC_+$,
 %$u_2(z):=u(z_1,z)$ is differentiable (in the complex sense) for all $z\in \CC_2$.
 }
 \end{defi}

 In the sequel we shall make use of the following simple,
 but useful, observation, saying that if two functions have
 identical Laplace transforms in an arbitrarily small interval
 then they coincide.
 \begin{lem}{\rm
 \label{lem.laplace}
 Let $I=(\lambda_1,\lambda_2)$ where $0\leq\lambda_1<\lambda_2\leq\infty$.
 Suppose that for the Borel functions $g_i:(0,\infty)\to \RR$, $i=1,2$, it is true
 that $\int_0^\infty e^{-\lambda x} |g_i(x)| \ud x<\infty$ for all $\lambda\in I$.
 If
 \[
 \int_0^\infty e^{-\lambda x} g_1(x) \ud x=\int_0^\infty e^{-\lambda x} g_2(x) \ud x,
 \ \ \ \lambda\in I,
 \]
 then $g_1(x)=g_2(x)$ for almost all $x\in(0,\infty)$.
 }
 \end{lem}
 \begin{pr}{Proof} For $i=1,2$ write
 $g_i=g_i^{+}-g_i^{-}$ where $g_i^{+}(x)=\max\{g_i(x),0\}$,
 $g_i^{-}(x)=\max\{-g_i(x),0\}$. By assumption,
 \[
 0\leq
 \int_0^\infty e^{-\lambda x} w_1(x) \ud x=\int_0^\infty e^{-\lambda x} w_2(x) \ud x<\infty,
 \ \ \ \lambda\in I,
 \]
 where $w_1=g_1^{+}+g_2^{-}\geq 0$, $w_2=g_1^{-}+g_2^{+}\geq 0$.
 Fix $\lambda_0\in I$
 and set $\nu_0:=\int_0^\infty e^{-\lambda_0 x} w_1(x) \ud x
 =\int_0^\infty e^{-\lambda_0 x} w_2(x) \ud x\geq 0$. Since $w_i\geq 0$, %$i=1,2$,
 it is clear that $\nu_0=0$ implies that
 $w_1=0$ and $w_2=0$ a.e.; then, $w_1=w_2$ a.e.\ and, thus,
 $g_1=g_2$ a.e.
 If $\nu_0>0$, we may define the probability densities
 $f_i(x)=e^{-\lambda_0 x}w_i(x)/\nu_0$, $i=1,2$; then, for $\lambda\in I$,
 \[
 \int_0^\infty e^{(\lambda_0-\lambda) x} f_1(x) \ud x
 =
 \frac{1}{\nu_0}\int_0^\infty e^{-\lambda x} w_1(x) \ud x
 =
 \frac{1}{\nu_0}\int_0^\infty e^{-\lambda x} w_2(x) \ud x
 =
 \int_0^\infty e^{(\lambda_0-\lambda) x} f_2(x) \ud x.
 %, \ \ \ \ \lambda\in I.
 %\int_0^\infty e^{-(\lambda-\lambda_0) x} \frac{e^{-\lambda_0 x} w_i(x)}{\nu_0} \ud x,
 %\ \ \ i=1,2.
 \]
 If $M_i$ is the moment generating function of $f_i$, it follows that
 $M_1(t)=M_2(t)<\infty$ for all $t\in J:=(\lambda_0-\lambda_2,\lambda_0-\lambda_1)$,
 and since $J$ contains the origin, $f_1=f_2$ a.e.; this implies that $w_1=w_2$ a.e.,
 and consequently, $g_1=g_2$ a.e.
 \medskip
 \end{pr}

 \begin{prop}{\rm
 \label{prop.continuity}
 Let $s(x,y):(0,\infty)^2\to\RR$ be a Borel measurable function
 satisfying
 \[
 \int_0^\infty \int_0^\infty
 |s(x,y)| \exp(-\lambda x-\alpha y) x^{n\alpha-1} \ud y \ud x<\infty
 \ \ \mbox{for all $\alpha>0$, $\lambda>0$.}
 \]
 Then, the function
 \[
 H(\alpha,\lambda):=
 \int_0^\infty \int_0^\infty
 s(x,y) \exp(-\lambda x-\alpha y) x^{n\alpha-1} \ud y \ud x
 \]
 is continuous in $(0,\infty)^2$.
 }
 \end{prop}
 \begin{pr}{Proof}
 Consider an arbitrary sequence
 $(\alpha_m,\lambda_m)\to (\alpha,\lambda)\in(0,\infty)^2$, as
 $m\to \infty$.
 Then, we can find $m_0$ such that
 $\frac{\alpha}{2}<\alpha_m<2\alpha$ and $\frac{\lambda}{2}<\lambda_m<2\lambda$
 for $m\geq m_0$,
 The function $\delta_1(\lambda)=e^{-\lambda x}$ is obviously decreasing
 (for each fixed $x>0$), so that $e^{-\lambda_m x}<e^{-\frac{\lambda}{2} x}$,
 $m\geq m_0$.
 Consider also the function $\delta_2(\alpha)=e^{-\alpha y} x^{n\alpha-1}$
 (for fixed $x>0$, $y>0$).
 Then, $\log \delta_2(\alpha)=(n\log x -y)\alpha -\log x$ is a linear
 function of $\alpha$, and thus,
 $
 %e^{-\alpha_m y} x^{n\alpha_m-1}=
 \delta_2(\alpha_m)<\max\{\delta_2(\alpha/2),\delta_2(2\alpha)\}<
 \delta_2(\alpha/2)+\delta_2(2\alpha)$, $m\geq m_0$.
 Combining the above we obtain the inequality
 \[
 e^{-\lambda_m x-\alpha_m y} x^{n\alpha_m-1}<
 e^{-\frac{\lambda}{2} x-\frac{\alpha}{2} y} x^{n\frac{\alpha}{2}-1}+
 e^{-\frac{\lambda}{2}x-2\alpha y} x^{2 n\alpha-1}, \ \ \ (x,y)\in(0,\infty)^2,
 \ \
 m\geq m_0.
 \]
 It follows that the sequence
 \[
 w_m(x,y):= s(x,y) e^{-\lambda_m x-\alpha_m y} x^{n\alpha_m-1},
 \ \ \ (x,y)\in(0,\infty)^2, \ \ m=m_0,m_0+1,\ldots,
 \]
 is dominated by $W(x,y):=|s(x,y)| e^{-\frac{\lambda}{2} x-\frac{\alpha}{2} y} x^{n\frac{\alpha}{2}-1}+ |s(x,y)|
 e^{-\frac{\lambda}{2}x-2\alpha y} x^{2 n\alpha-1}$,
 and, by assumption, $W$ is integrable on $(0,\infty)^2$.
 Clearly, $w_m(x,y)\to w(x,y):=
 s(x,y) e^{-\lambda x-\alpha y} x^{n\alpha-1}$ as $m\to\infty$.
 Thus, by dominated convergence,
 \[
 H(\alpha_m,\lambda_m)=\int_0^\infty \int_0^\infty
 w_m(x,y) \ud y \ud x \to
 \int_0^\infty \int_0^\infty
 w(x,y) \ud y \ud x =H(\alpha,\lambda), \ \ \ \mbox{as} \ \ m\to\infty,
 \]
 completing the proof.
 \end{pr}

 \begin{prop}{\rm
 \label{prop.h(a)}
 Let $h(\alpha):(0,\infty)\to\RR$ be a parametric function
 of the shape parameter. If there exists an unbiased estimator $U_1=U_1(\bbb{X})$
 for $h$, then the UMVUE of $h(\alpha)$ is a function of $Y$, alone.
 }
 \end{prop}
 \begin{pr}{Proof}
 By assumption, $\E_{(\alpha,\lambda)} |U_1(\bbb X)|<\infty$ for all
 $\alpha>0$, $\lambda>0$.
 The function $U_1:(0,\infty)^n\to \RR$ can always be taken
 to be Borel measurable. From the RB/LS Theorem,
 the UMVUE is given by $U(X,Y)=\E(U_1|X,Y)$,
 where $U:(0,\infty)^2\to\RR$ can also be taken to be Borel measurable
 with no loss of generality. Unbiasedness of the estimator
 $U(X,Y)$ means that (recall that $X,Y$ are independent)
 \[
 \int_0^{\infty} \int_{0}^1 \left\{\frac{\lambda^{n\alpha}}{\Gamma(n\alpha)}
 x^{n\alpha-1} e^{-\lambda x}\right\}\left\{
 K(\alpha) y^{\alpha-1}(1-y)^{(n-3)/2} G(y) \right\}
 U(x,y) \ud y
 \ud x = h(\alpha),
 \]
 for all $\alpha>0$, $\lambda>0$. This relation can be rewritten in the form
 \[
 H(\alpha,\lambda):=\frac{h(\alpha) \Gamma(n\alpha)}{K(\alpha)\lambda^{n\alpha}}=
 \int_0^{\infty} \int_{0}^\infty (1-e^{-y})^{(n-3)/2} G(e^{-y}) U(x,e^{-y})
  \exp(-\lambda x-\alpha y) x^{n\alpha-1}
  \ud y
 \ud x.
 \]
 Since $\E_{(\alpha,\lambda)} |U(X,Y)|<\infty$ for all $\alpha>0$ and $\lambda>0$, we
 may apply Proposition \ref{prop.continuity} with
 $s(x,y)=(1-e^{-y})^{(n-3)/2} G(e^{-y}) U(x,e^{-y})$
 to conclude that the function
 $H(\alpha,\lambda)$
 is continuous, and hence, $h(\alpha)$ is a continuous function
 of $\alpha>0$; notice that all functions $K(\alpha)$, $\Gamma(n\alpha)$,
 and $\lambda^{n\alpha}$, are positive and continuous.

 Write now the relation $h(\alpha)=\E_{(\alpha,\lambda)} U(X,Y)$
 in the form
 \[
 %\label{laplace}
 \int_0^\infty e^{-\lambda x} h(\alpha) x^{n \alpha-1}\ud x =
 \int_0^{\infty} e^{-\lambda x} K(\alpha) x^{n\alpha-1} \left\{
 \int_{0}^1 y^{\alpha-1}(1-y)^{(n-3)/2} G(y) U(x,y) \ud y
 \right\} \ud x.
 \]
 %and fix $\alpha=\alpha_0>0$.
 %Then, (\ref{laplace})
 %, which holds for  all $\lambda>0$,
 %shows
 It follows that the Laplace
 transforms of the functions $f_1(x):= h(\alpha) x^{n \alpha-1}$
 and $f_2(x):= K(\alpha) x^{n\alpha-1}
 \int_{0}^1 y^{\alpha-1}(1-y)^{(n-3)/2} G(y) U(x,y) \ud y$
 are identical,
 %From the one-to-one property of the Laplace transform
 and Lemma \ref{lem.laplace} implies that for every $\alpha>0$,
 \be
 \label{a.e.}
 \frac{h(\alpha)}{K(\alpha)}=\int_{0}^1 y^{\alpha-1}(1-y)^{(n-3)/2} G(y)
 U(x,y)
 \ud y, \ \ \  \mbox{for almost all} \ x\in(0,\infty).
 \ee
 Write $\QQ_+=\QQ\cap(0,\infty)$
 and set $A_\alpha=\{x>0: \mbox{(\ref{a.e.}) holds}\}$,
 $N_\alpha=(0,\infty)\setminus A_\alpha=\{x>0: \mbox{(\ref{a.e.}) does not hold}\}$,
 %$B=(0,\infty)\cap \QQ$,
 $A=\cap_{\alpha\in \QQ_+}A_\alpha$,
 $N=\cup_{\alpha\in \QQ_+}N_\alpha=(0,\infty)\setminus A$.
 By construction, the Lebesgue  measure of $N$ is zero.
 %$|N|=0$, where $|N|$ denotes the Lebesgue measure of $N$.
 We chose an arbitrary $x_0\in A$ and define $\widetilde{u}(y):=U(x_0,y)$.
 Notice that $\widetilde{u}:(0,\infty)\to\RR$ is not necessarily
 Borel measurable;
 it is, however, Lebesgue measurable and, hence, it is equal
 a.e.\ to a Borel function $u:(0,\infty)\to\RR$. Obviously,
 substituting $u$ in place of $\widetilde{u}$, the values
 of the above integrals are not affected.
 %since $U:(0,\infty)^2\to\RR$ is Borel.
 From (\ref{a.e.})
 we have
 \be
 \label{contin}
 \frac{h(\alpha)}{K(\alpha)}=\int_{0}^1 y^{\alpha-1}(1-y)^{(n-3)/2} G(y) u(y)
 \ud y, \ \ \ \ \alpha\in \QQ_+.
 \ee
 The right-hand side of (\ref{contin}) defines a continuous function
 of $\alpha>0$, because
 \[
 \int_{0}^1 y^{\alpha-1}(1-y)^{(n-3)/2} G(y) u(y) \ud y
 =
 \int_{0}^{\infty} e^{-\alpha y} (1-e^{-y})^{(n-3)/2}
 G(e^{-y}) u(e^{-y}) \ud y
 \]
 is the Laplace transform of the function
 $y\mapsto (1-e^{-y})^{(n-3)/2} G(e^{-y}) u(e^{-y})$.
 Since $h(\alpha)$ and $K(\alpha)$ are both continuous,
 it follows from (\ref{contin}) that
 $h(\alpha)=
 \int_{0}^1 K(\alpha) y^{\alpha-1}(1-y)^{(n-3)/2} G(y) u(y) \ud y=
 \int_0^1 f_Y(y;\alpha)u(y) \ud y$
 {\it for all} $\alpha>0$. Hence, $u(Y)$ is an unbiased
 estimator of $h(\alpha)$, and also, it is a function
 of the complete, sufficient statistic $(X,Y)$; thus,
 it is the (unique) UMVUE.
 \end{pr}

 \begin{rem}{\rm
 \label{rem.FI} (A consequence to the Cram\'{e}r-Rao bound).
 The Fisher information matrix based on a single observation
 from $\mathcal{G}(\alpha,\lambda)$ is given by
 \[
 J=J(\alpha,\lambda)=\frac{1}{\lambda^2}\left(
 \begin{array}{cc}
 \lambda^2 \psi_1(\alpha) & -\lambda \\
 -\lambda & \alpha
 \end{array}
 \right),
 \mbox{ with inverse }
 J^{-1}=\frac{1}{\alpha \psi_1(\alpha)-1}\left(
 \begin{array}{cc}
 \alpha & \lambda \\
 \lambda & \lambda^2 \psi_1(\alpha)
 \end{array}
 \right),
 \]
 where $\psi_1(\alpha)=(\log \Gamma(\alpha))''>1/\alpha$.
 The Cram\'{e}r-Rao (CR) bound says that for any unbiased estimator
 $T=T(\bbb X)$ of the parametric function $h(\alpha,\lambda)$,
 $n \Var T\geq (\nabla h(\alpha,\lambda))^T J^{-1}\nabla h(\alpha,\lambda)$.
 Hence, in the particular case where $h=h(\alpha)$ is a parametric function
 of the shape parameter alone
 (i.e., $h$ does not depend on $\lambda$),
 the CR bound takes
 the form
 \be
 \label{CR_alpha_lambda}
 n\Var T(\bbb X) \geq \frac{h'(\alpha)^2}{\psi_1(\alpha)-1/\alpha},
 \ \ \ \mbox{whenever } \ \E_{(\alpha,\lambda)} T(\bbb X)=h(\alpha)
 \ \mbox{ for all } \alpha>0, \ \lambda>0.
 \ee
 We shall now show that Proposition \ref{prop.h(a)} yields a
 better lower bound. Indeed, since the UMVUE of $h(\alpha)$
 is a Borel function of $Y$, say $u=u(Y)$,
 we can apply the CR bound, $\Var u(Y)\geq h'(\alpha)^2/J_Y$,
 where $J_Y$ is the Fisher information of $Y$.
 From (\ref{f_Y}), the density of $Y$ is given by
 \[
 f_Y(y;\alpha)=
 \exp\Big(a\log y-A(\alpha)\Big)k(y),
 \]
 where $A(\alpha)=-\log K(\alpha)$ and
 $k(y)=y^{-1}(1-y)^{(n-3)/2}G(y) I_{(0,1)}(y)$,
 with $I_B$ denoting the indicator function of $B$. This shows that
 $\{ f_Y(y;\alpha), \alpha>0\}$ defines
 a natural exponential family of Lebesgue densities, and thus,
 the regularity conditions
 are fulfilled.
 We calculate $-\frac{\partial^2}{\partial\alpha^2}\log f_Y(y;\alpha)=A''(\alpha)$
 and, therefore,
 \begin{eqnarray*}
 J_Y
 \hspace*{-1ex}& = & \hspace*{-1ex}
 (-\log K(\alpha))''=\left(n \log \Gamma(\alpha)
 +\sum_{i=1}^{n-1}\log \Gamma(i/n)-\sum_{i=1}^{n-1}\log \Gamma(\alpha+i/n)\right)''
 \\
 \hspace*{-1ex}& = & \hspace*{-1ex}
 n\psi_1(\alpha)-\sum_{i=1}^{n-1}\psi_1(\alpha+i/n).
 \end{eqnarray*}
 It follows that for any unbiased estimator $T(\bbb X)$ of $h(\alpha)$,
 \be
 \label{CR_alpha}
 n \Var T(\bbb X) \geq n \Var u(Y)\geq \frac{h'(\alpha)^2}{\psi_1(\alpha)-n^{-1}\sum_{i=1}^{n-1}\psi_1(\alpha+i/n)}.
 \ee
 Since the function $\psi_1(\alpha)$ is positive (the $\Gamma$ function
 is log-convex) and decreasing, we have
 $n^{-1}\psi_1(\alpha+i/n)<\int_{(i-1)/n}^{i/n}\psi_1(\alpha+t)\ud t$,
 $i=1,\ldots,n$. Therefore,
 \begin{eqnarray*}
 \frac{1}{n} \sum_{i=1}^{n-1}\psi_1(\alpha+i/n)<
 \frac{1}{n} \sum_{i=1}^{n}\psi_1(\alpha+i/n)<
 \int_0^1 \psi_1(\alpha+t) \ud t
 =\int_\alpha^{\alpha+1} \psi_1(t) \ud t
 \\
 =(\log \Gamma(\alpha+1))'-(\log \Gamma(\alpha))'
 =\left(\log \frac{\Gamma(\alpha+1)}{\Gamma(\alpha)}\right)'
 =(\log \alpha)'=\frac{1}{\alpha}.
 \end{eqnarray*}
 Thus, excluding the trivial case where
 $h(\alpha)$ is the constant function, the lower bound in (\ref{CR_alpha})
 is strictly larger than the CR bound of (\ref{CR_alpha_lambda}).
 This means that,
 for any given non-constant parametric function $h(\alpha)$,
 no efficient estimator exists, when efficiency is defined
 in the traditional way, i.e., in terms of attainability
 of the bound in (\ref{CR_alpha_lambda}). In the contrary,
 there are functions of the shape parameter that can be efficiently
 estimated in the sense of (\ref{CR_alpha}). More precisely,
 it can be checked that for any constants $c_1,c_2$,
 the parametric function
 $h(\alpha)=c_1\Big(\psi(\alpha)-n^{-1}\sum_{i=1}^{n-1}\psi(\alpha+i/n)\Big)+c_2$,
 where $\psi(\alpha)=\Gamma'(\alpha)/\Gamma(\alpha)$,
 is efficiently estimated by $T(\bbb X)=u(Y)=(c_1/n)\log Y +c_2$,
 in the sense that $\E_{(\alpha,\lambda)} T(\bbb X)=h(\alpha)$
 for all $\alpha,\lambda$, and
 $n \Var T(\bbb X)=c_1^2 \Big(\psi_1(\alpha)-
 n^{-1}\sum_{i=1}^{n-1}\psi_1(\alpha+i/n)\Big)$,
 which is equal to the lower bound given by (\ref{CR_alpha}).
 Finally, it can be easily verified that all efficiently estimated
 functions are of the above form. Notice, however, that
 asymptotically the lower bounds do coincide, since,
 by the Riemann integral,
  $\lim_{n\to\infty} n^{-1} \sum_{i=1}^{n-1}\psi_1(\alpha+i/n)
 =\int_\alpha^{\alpha+1}\psi_1(t) \ud t=1/\alpha$.
 }
 \end{rem}

 We now state the main result of this section. To the best
 of our knowledge, this result, as well as the existence
 of $u_0$, below, is new.

 \begin{theo}{\rm
 \label{theo.unbiased}
 Let $X_1,\ldots,X_n$, $n\geq 2$, be a random sample
 from $\mathcal{G}(\alpha,\lambda)$.
 For every $n\geq 4$, the UMVUE of the shape parameter $\alpha$ is
 given by
 \be
 \label{u_0}
 u_0(Y):= \frac{1}{2}(n-3)\frac{ Y}{1-Y}-\frac{Y G'(Y)}{G(Y)},
 \ee
 where $G$ is as in (\ref{g_(n-1)}) and $Y^{1/n}$
 is the ratio of the geometric to the arithmetic mean
 of $X_1,\ldots,X_n$. The estimator $u_0(Y)$
 has finite variance if and only if $n\geq 6$.
 Moreover, for $n=2,3$, no unbiased estimator of $\alpha$ exists.
 }
 \end{theo}
 \begin{pr}{Proof}
 According to Proposition \ref{prop.h(a)}, if an unbiased estimator
 (of $\alpha$) exists, the UMVUE must be a function of $Y$, say $u(Y)$.
 Then, the relation $\E_{\alpha} u(Y)=\alpha$ can be written as
 \[
 \frac{1}{K(\alpha)}=\frac{1}{\alpha}
 \int_0^1 y^{\alpha-1} (1-y)^{(n-3)/2} G(y) u(y) \ud y,
 \ \ \alpha>0.
 \]
 Observing that $1/K(\alpha)=\int_0^1 y^{\alpha-1} (1-y)^{(n-3)/2} G(y) \ud y$ and
 $1/\alpha=\int_0^{\infty}e^{-\alpha x} \ud x$,
 the substitution $y=e^{-x}$ yields
 \[
 \int_0^\infty e^{-\alpha x} (1-e^{-x})^{(n-3)/2} G(e^{-x}) \ud x =
 \left\{\int_0^{\infty}e^{-\alpha x} \ud x\right\}
 \left\{\int_0^\infty e^{-\alpha x} (1-e^{-x})^{(n-3)/2}
 G(e^{-x}) u(e^{-x}) \ud x\right\},
 \]
 for all $\alpha>0$. Since we have assumed that
 $\E_{\alpha}|u(Y)|<\infty$ for all $\alpha>0$, we may apply Fubini's theorem
 to the right-hand side of the above equation, obtaining
 \[
 \int_0^\infty e^{-\alpha x} (1-e^{-x})^{(n-3)/2} G(e^{-x}) \ud x =
 \int_0^{\infty}e^{-\alpha x} \left\{\int_0^x  (1-e^{-t})^{(n-3)/2}
 G(e^{-t}) u(e^{-t}) \ud t\right\} \ud x,
 \ \ \ \alpha>0.
 \]
 This relation shows that the Laplace transforms of
 $f_1(x):=(1-e^{-x})^{(n-3)/2} G(e^{-x})$ and
 $f_2(x):=\int_0^x  (1-e^{-t})^{(n-3)/2}
 G(e^{-t}) u(e^{-t}) \ud t$ are identical.
 Hence, $f_1\equiv f_2$ (notice that both functions are continuous).
 Setting $x=-\log y$ we get
 \be
 \label{int.equation}
 (1-y)^{(n-3)/2} G(y)=\int_0^{-\log y}  (1-e^{-t})^{(n-3)/2}
 G(e^{-t}) u(e^{-t}) \ud t = \int_{y}^{1}
 \frac{1}{x}  (1-x)^{(n-3)/2}
 G(x) u(x) \ud x,
 \ee
 $0<y<1$, and a differentiation of (\ref{int.equation}) leads to
 \[
  -\frac{n-3}{2} (1-y)^{(n-5)/2} G(y) +(1-y)^{(n-3)/2} G'(y)=
  \frac{-u(y)}{y} (1-y)^{(n-3)/2}
  G(y), \ \ \ \mbox{for almost all}\ y\in(0,1).
 \]
 Hence, solving for $u$ we see that $u(y)=\big((n-3)/2\big)y/(1-y)-yG'(y)/G(y)=u_0(y)$,
 a.e., with $u_0$ as in (\ref{u_0}). The preceding argument shows that
 if an unbiased estimator $T(X_1,\ldots,X_n)$ exists then $u_0(Y)$ must be unbiased,
 but
 {\it it does not prove
 that an unbiased estimator exists}. To see this, let $n=2$.
 Then, since $G(y)\equiv 1$ when $n=2$,
 (\ref{int.equation}) reads as
 \be
 \label{count_n=2}
 \frac{1}{\sqrt{1-y}}=\int_{y}^1 \frac{u(x)}{x\sqrt{1-x}} \ud x.
 \ee
 The assumption $\E_{\alpha} |u(Y)|<\infty$ is equivalent
 to $\int_0^1 x^{\alpha-1} (1-x)^{-1/2}|u(x)|\ud x<\infty$ for all $\alpha>0$.
 Hence,
 \[
 \left|\int_{y}^1 \frac{u(x)}{x\sqrt{1-x}} \ud x \right| \leq
 \int_{y}^1 \frac{|u(x)|}{x\sqrt{1-x}} \ud x =
 \int_{y}^1 \frac{x^{\alpha-1}|u(x)|}{\sqrt{1-x}} \frac{1}{x^\alpha} \ud x
 \leq \frac{1}{y^\alpha}\int_{y}^1 \frac{x^{\alpha-1}|u(x)|}{\sqrt{1-x}} \ud x \to 0,
 \]
 \mbox{as} $y\to 1-$, and this contradicts (\ref{count_n=2}).
 When $n=3$, it is easy to see that the function $G$ of (\ref{g_(n-1)})
 is given by
 \[
 G(y)=\sum_{m=0}^\infty \frac{(3 m)!}{(3^m m!)^3} (1-y)^m, \ \ \ 0<y<1.
 \]
 Hence, $G(1-)=1$. Since $(1-y)^{(n-3)/2}\equiv 1$, assuming
 $\E_{\alpha} |u(Y)|<\infty$ we get
 $\int_0^1 x^{\alpha-1} G(x) |u(x)|\ud x<\infty$,
 $\alpha>0$. Therefore,
 \[
 \left|\int_{y}^1 \frac{G(x)}{x}u(x) \ud x \right| \leq
 \int_{y}^1 \frac{G(x)}{x}|u(x)| \ud x =
 \int_{y}^1 \frac{x^{\alpha-1}G(x)|u(x)|}{x^\alpha} \ud x
 \leq \frac{1}{y^{\alpha}}
 \int_{y}^1 x^{\alpha-1} G(x) |u(x)| \ud x \to 0,
 \]
 \mbox{as} $y\to 1-$.
 This contradicts (\ref{int.equation}), i.e.,
 \[
 G(y)=\int_{y}^1 \frac{1}{x}G(x)u(x) \ud x,
 \]
 since the left-hand side of the above equation approaches $1$ as $y\to 1-$.

 The preceding argument shows that there is no unbiased estimator when
 $n\leq 3$. For $n\geq 4$, however, the situation is completely
 different: The (positive)
 estimator $u_0(Y)$ of (\ref{u_0}) has finite expectation for all $\alpha>0$ and,
 indeed, $\E_{\alpha} u_0(Y)=\alpha$. We now proceed to verify this claim.
 First observe that
 $-G'(y)>0$ so that $u_0(y)>0$ for $0<y<1$. Noting that
 $G'(y)=-\sum_{m=0}^\infty m d_m (1-y)^{m-1}$,
 we calculate
 \begin{eqnarray*}
 \E_{\alpha} u_0(Y)
 &\hspace*{-1ex}=\hspace*{-1ex}&
 K(\alpha) \int_0^1 y^{\alpha} (1-y)^{(n-5)/2}
 \left(\sum_{m=0}^\infty \left(\frac{n-3}{2}+m\right) d_m (1-y)^m\right) \ud y
 \\
 &\hspace*{-1ex}=\hspace*{-1ex}&
 K(\alpha) \sum_{m=0}^\infty \left(\frac{n-3}{2}+m\right)d_m
 \int_0^1 y^\alpha (1-y)^{m+(n-5)/2} \ud y
 \\
 &\hspace*{-1ex}=\hspace*{-1ex}&
 \alpha\Gamma(\alpha) K(\alpha)
 \sum_{m=0}^\infty \frac{\Gamma(m+(n-1)/2)}{\Gamma(\alpha+m+(n-1)/2)} d_m,
 \end{eqnarray*}
 where, since $d_m>0$, the
 interchanging of summation and integration is justified by
 Beppo Levi's theorem. It remains to verify that
 \be
 \label{unbiasd}
 S:=\sum_{m=0}^\infty \frac{\Gamma(m+(n-1)/2)}{\Gamma(\alpha+m+(n-1)/2)} d_m
 =\frac{1}{\Gamma(\alpha) K(\alpha)}, \ \ \ \alpha>0.
 \ee
 From the definition of $d_m=\gamma_{n-1}(m)$ we have (set $i=n-2$ in (\ref{gamma}))
 \[
 d_m= \gamma_{n-1}(m)=
 \frac{\Gamma(m+(n-1)(n-2)/(2n))}
 {\Gamma(m+(n-1)/2)}
 \sum_{k=0}^m \frac{\Gamma(m-k+(n-1)/n)}{(m-k)!} \gamma_{n-2}(k).
 \]
 Substituting this expression to the series in (\ref{unbiasd})
 and changing the order of summation according to Tonelli's theorem,
 we obtain
 \begin{eqnarray}
 \nonumber
 S
 &\hspace*{-1ex}=\hspace*{-1ex}&
 \sum_{k=0}^\infty   \gamma_{n-2}(k)
 \sum_{m=k}^\infty
 \frac{\Gamma(m+(n-1)(n-2)/(2n))\ \Gamma(m-k+(n-1)/n)}
 {(m-k)!\ \Gamma(\alpha+m+(n-1)/2)}
 \\
 \nonumber
 &\hspace*{-1ex}=\hspace*{-1ex}&
 \sum_{k=0}^\infty   \gamma_{n-2}(k)
 \sum_{m=0}^\infty
 \frac{\Gamma(m+k+(n-1)(n-2)/(2n))\ \Gamma(m+(n-1)/n)}
 {m!\ \Gamma(m+\alpha+k+(n-1)/2)}
 \\
 \label{S}
 &\hspace*{-1ex}=\hspace*{-1ex}&
 \frac{\Gamma(\alpha)\ \Gamma((n-1)/n)}{\Gamma(\alpha+(n-1)/n)}
 \sum_{k=0}^\infty
 \frac{\Gamma(k+(n-1)(n-2)/(2n))}{\Gamma(\alpha+k+(n-1)(n-2)/(2n))} \gamma_{n-2}(k),
 \end{eqnarray}
 where for the last equality we applied the identity
 \[
 \sum_{m=0}^\infty \frac{\Gamma(m+\rho_1)\ \Gamma(m+\rho_2)}{m!\
 \Gamma(m+\rho+\rho_1+\rho_2)}=\frac{\Gamma(\rho)\ \Gamma(\rho_1)\ \Gamma(\rho_2)}
 {\Gamma(\rho+\rho_1) \ \Gamma(\rho+\rho_2)}, \ \ \ \ \ \rho,\ \rho_1,\ \rho_2>0,
 \]
 with $\rho=\alpha$, $\rho_1=k+(n-1)(n-2)/(2n)$, $\rho_2=(n-1)/n$.
 Since the function $f_{Y_1}(y;\alpha)$ in (\ref{f_Y1}) is a density, we have
 (see (\ref{g1}))
 \[
 \frac{1}{K_1(\alpha)}=\int_{0}^1 y^{\alpha-1}(1-y)^{(n-1)(n-2)/(2n)-1} G_1(y) \ud y
 =\sum_{k=0}^\infty \gamma_{n-2}(k)
 \int_{0}^1 y^{\alpha-1}(1-y)^{(n-1)(n-2)/(2n)+k-1} \ud y,
 \]
 and hence,
 \[
 \sum_{k=0}^\infty
 \frac{\Gamma(k+(n-1)(n-2)/(2n))}{\Gamma(\alpha+k+(n-1)(n-2)/(2n))} \gamma_{n-2}(k) =\frac{1}{\Gamma(\alpha)K_1(\alpha)}.
 \]
 Using the fact that
 $K_1(\alpha)=K(\alpha) \Gamma(\alpha) \Gamma((n-1)/n)/\Gamma(\alpha+(n-1)/n)$ and substituting the above equality in (\ref{S}) we conclude that
 %that $K_1(\alpha)=\Gamma(\alpha)^{-(n-1)}\prod_{i=1}^{n-2}
 %\Gamma(\alpha+i/n)/\Gamma(i/n)$, $K(\alpha)=\Gamma(\alpha)^{-n}\prod_{i=1}^{n-1}
 %\Gamma(\alpha+i/n)/\Gamma(i/n)$, we obtain
 %
 %$K_1(\alpha)/K(\alpha)=\Gamma(\alpha)^{-(n-1)}\left(\prod_{i=1}^{n-2}
 %\Gamma(\alpha+i/n)/\Gamma(i/n)\right) \Gamma(\alpha)^{n}\left(\prod_{i=1}^{n-1}
 %\Gamma(i/n)/\Gamma(\alpha+i/n)\right)$
  \[
 S=
 \frac{\Gamma((n-1)/n)}{\Gamma(\alpha+(n-1)/n)\ K_1(\alpha)}=
 \frac{\Gamma((n-1)/n) \ \Gamma(\alpha+(n-1)/n)}{
 K(\alpha)\ \Gamma(\alpha) \ \Gamma((n-1)/n)\ \Gamma(\alpha+(n-1)/n)}
 =\frac{1}{K(\alpha)\Gamma(\alpha)};
 \]
 this verifies (\ref{unbiasd}) and the claim is proved.
 Thus, for $n\geq 4$,
 $u_0(Y)$ is the UMVUE of the parametric function $h(\alpha)=\alpha$.

 We shall now investigate the finiteness of the variance of the UMVUE.
 Assume first that $n\geq 6$.
 From the inequality
 %It remains to show that $\E_{\alpha} u_0(Y)^2<\infty$ if and only if
 %$n\geq 6$. Since
 $u_0(y)^2\leq 2c^2 y^2/(1-y)^2+2 y^2 |G'(y)|^2/G(y)^2$, where $c=(n-3)/2$,
 %it suffices to show that
 we see that $\E_{\alpha} u_0(Y)^2<\infty$ if both integrals
 \[
 \int_{0}^1 y^{\alpha+1} (1-y)^{(n-7)/2} G(y) \ud y
 \ \ \ \mbox{and}
 \ \ \
 \int_{0}^1 y^{\alpha+1} (1-y)^{(n-3)/2} \ \frac{|G'(y)|^2}{G(y)} \ud y
 \]
 are finite. The integrands are continuous functions in $(0,1)$, hence,
 it suffices to verify their integrability in a neighborhood
 of the endpoints. Obviously,
 \[
 \int_{1-\epsilon}^1 y^{\alpha+1} (1-y)^{(n-7)/2} G(y) \ud y<\infty
 \ \ \ \mbox{and}
 \ \ \
 \int_{1-\epsilon}^1 y^{\alpha+1} (1-y)^{(n-3)/2} \ \frac{|G'(y)|^2}{G(y)}
 \ud y <\infty,
 \]
 since $n\geq 6$, $\lim_{y\to 1-} y^{\alpha+1} G(y)=d_0$
 and $\lim_{y\to 1-}y^{\alpha+1} |G'(y)|^2/G(y)=d_1^2/d_0$,
 so that, close to $y=1$, the first integrand behaves like $(1-y)^{(n-7)/2}$
 and the second one like $(1-y)^{(n-3)/2}$.
 Regarding the behavior of the first integral close to $y=0$
 it suffices to observe that
 %we proceed as follows: Since
 $G(y)$ is positive
 %decreasing,
 and $\int_{0}^1 y^{\alpha-1} (1-y)^{(n-3)/2} G(y)\ud y=1/K(\alpha)<\infty$;
 see (\ref{f_Y}).
 Then, we trivially have
 %$\lim_{y\to 0+}
 $\int_{0}^\delta y^{\alpha+1} (1-y)^{(n-7)/2} G(y)\ud y
 %< \int_{0}^\epsilon y^{\alpha-1} (1-y)^{(n-3)/2} (1-y)^{-2} g(y)\ud y
 <(1-\delta)^{-2} \int_{0}^\delta y^{\alpha-1} (1-y)^{(n-3)/2} G(y)\ud y<\infty$.
 %(1-\epsilon)^{-(n-3)/2}
 %\int_{0}^\epsilon y^{\alpha-1} (1-y)^{(n-3)/2} g(y)\ud y<\infty$.
 %\leq \lim_{y\to 0+}\int_{0}^y x^{\alpha-1} (1-x)^{(n-3)/2} g(x)\ud x=0$,
 %and thus,
 %\[
 %y^{\alpha} g(y) =\frac{\alpha}{2^\alpha-1} g(y)\int_{y}^{2y} x^{\alpha-1} \ud x
 %< \frac{\alpha}{2^\alpha-1} \int_{y}^{2y} x^{\alpha-1} g(x) \ud x \to 0,
 %\ \ \ \mbox{ as } \ \ y\to 0+.
 %\]
 It remains to verify that
 $\int_0^{\delta} y^{\alpha+1} (1-y)^{(n-3)/2} \
 |G'(y)|^2/G(y) \ud y <\infty$, and since $1/G(y)<1/d_0$ and
 $(1-y)^{(n-3)/2}<1$,
 it suffices to show that for some $\delta\in(0,1)$,
 \be
 \label{var}
 \int_0^{\delta} y^{\alpha+1}  |G'(y)|^2 \ud y <\infty.
 \ee
 %To this end,
 Recall that $\E_{\alpha} u_0(Y)<\infty$; since
 $u_0(y)>y |G'(y)|/G(y)$, it follows that
 \[
 \int_0^1 y^{\alpha}(1-y)^{(n-3)/2} |G'(y)| \ud y<\infty,
 \]
 and thus,
 $\int_0^{1/2} y^{\alpha} |G'(y)| \ud y<\infty$.
 The function $|G'(y)|=\sum_{m=1}^\infty m d_{m} (1-y)^{m-1}$
 is strictly decreasing in $(0,1)$
 and therefore, with $c_\alpha=(\alpha+1)/(1-2^{-\alpha-1})>0$,
 we get
 %; since,
 %by dominated convergence,
 %$\int_{y/2}^{y} x^{\alpha}|g'(x)| \ud x<\int_0^{y} x^{\alpha}|g'(x)| \ud x
 %\to 0$ as $y\to 0+$, we obtain
 \[
 y^{\alpha+1}|G'(y)|=c_\alpha \ |G'(y)|
 \int_{y/2}^{y} x^\alpha \ud x
 <c_\alpha  \int_{y/2}^{y} x^\alpha |G'(x)| \ud x
 <c_\alpha \int_{0}^{y} x^\alpha |G'(x)| \ud x \to 0,
 \]
 as $y\to 0+$, by dominated convergence.
 Hence, $\lim_{y\to 0+}y^{\alpha+1}|G'(y)|=0$ (for all $\alpha>0$).
 Thus, for every $\epsilon>0$ we can find $\delta=\delta(\alpha,\epsilon)>0$
 such that $\big(y^{\alpha/4+1}|G'(y)|\big)^2<\epsilon^2$ for
 $y\in(0,\delta)$. Then, for this value of $\delta$,
 \[
 \int_{0}^\delta y^{\alpha+1}  |G'(y)|^2 \ud y
 =
 \int_{0}^\delta y^{\alpha/2-1}  \Big(y^{\alpha/4+1}|G'(y)|\Big)^2 \ud y
 <
 \epsilon^2 \int_{0}^\delta y^{\alpha/2-1}  \ud y<\infty,
 \]
 and (\ref{var}) is proved. We thus conclude that the variance of $u_0(Y)$
 is finite for all $n\geq 6$ and $\alpha>0$.

 Finally, assume that $n=4$ or $5$.
 Since $u_0(y)=c y/(1-y)+y|G'(y)|/G(y)>c y/(1-y)$, $c=(n-3)/2>0$,
 and $G(y)>d_0>0$,  the variance of the UMVUE is infinite:
 \[
 \E_{\alpha} u_0(Y)^2 >
 c^2 K(\alpha) \int_{0}^1 y^{\alpha+1} (1-y)^{(n-7)/2} G(y)  \ud y
 >  c^2 K(\alpha) d_0 \int_{0}^1 y^{\alpha+1} (1-y)^{(n-7)/2} \ud y=\infty
 \]
 for all $\alpha>0$; this completes the proof of the theorem.
 \end{pr}

 \begin{cor}{\rm
 \label{cor.lambda} For $n\geq 4$, the UMVUE of the rate parameter (reciprocal scale parameter) $\lambda$ is given by
 \[
 u_1(X,Y)=\frac{u_0(Y)-1/n}{\overline{X}}
 = \frac{n \ u_0(Y)-1}{X}, \ \ \ \ n>1/\alpha,
 \]
 where $u_0(Y)$ is as in (\ref{u_0}) and
 $\overline{X}=X/n=n^{-1}\sum_{i=1}^n X_i$.
 For $n\geq 6$, the estimator $u_1(X,Y)$
 has finite variance if and only if $\alpha>2/n$ and $\lambda>0$,
 while for $n=4$ or $5$ its variance is infinite. Finally,
 if either $n\in\{2,3\}$ or $n\geq 4$ and $0<\alpha\leq n^{-1}$,
 no unbiased estimator
 for $\lambda$ exists.
 }
 \end{cor}

 \begin{pr}{Proof}
  % of Corollary \ref{cor.lambda}}
 By definition, $u_1$ is a function of the complete, sufficient
 statistic $(X,Y)$. Since $X,Y$ are independent and
 $\E_{(\alpha,\lambda)} X^{-1} = \lambda/(n \alpha-1)$, $n>\alpha^{-1}$,
 we obtain
 $\E_{(\alpha,\lambda)} u_1(X,Y)=\Big(n\E_{(\alpha,\lambda)}u_0(Y)-1\Big)
 \left(\E_{(\alpha,\lambda)} X^{-1}\right)  =\lambda$,
 %Moreover, the function $u_1(z_1,z_2) =(n \ u_0(z_2)-1)/z_1$
 %is, obviously, holomorphic in $\CC_+\times \CC_1$, and
 %the set $\{(\alpha,\lambda):\alpha>1/n, \lambda>0\}$
 %has nonempty interior, showing that
 and $u_1(X,Y)$ is, indeed, the UMVUE, according to Definition
 \ref{def.UMVUE}. Let $n=4$ or $5$.
 Since $u_1(x,y)^2=x^{-2}(n^2 u_0(y)^2-2nu_0(y)+1)$
 and $\E_{(\alpha,\lambda)}u_0(Y)^2=\infty$, it follows
 that
 %$2n \alpha+n^2\E_{(\alpha,\lambda)}u_0(Y)^2=
 $\E_{(\alpha,\lambda)}\Big(n^2 u_0(Y)^2-2nu_0(Y)+1\Big)=\infty$
 and, thus, by the independence of $X,Y$,
 $\E_{(\alpha,\lambda)} u_1(X,Y)^2
 =\left(\E_{(\alpha,\lambda)} X^{-2}\right) \left(\E_{(\alpha,\lambda)}\Big(n^2 u_0(Y)^2-2nu_0(Y)+1\Big)\right)=\infty$.
 Next, assume $n\geq 6$. In this case,
 $\E_{(\alpha,\lambda)} u_1(X,Y)^2
 \leq \left(\E_{(\alpha,\lambda)} X^{-2}\right)
 \left(1+n^2\E_{(\alpha,\lambda)}u_0(Y)^2\right)$,
 and since $\E_{(\alpha,\lambda)}u_0(Y)^2<\infty$, we obtain
 $\E_{(\alpha,\lambda)} u_1(X,Y)^2<\infty$ if and only if
 $\E_{(\alpha,\lambda)} X^{-2}<\infty$; equivalently, $\alpha>2 n^{-1}$.

 Assume now that $n\geq 2$, $0<\alpha\leq n^{-1}$,
 and set $\theta:=n \alpha\in(0,1]$
 for fixed $\alpha$. If for this particular value
 of $\alpha$ the estimator $w(X,Y)$
 were unbiased for $\lambda>0$, then we could find
 an interval
 $I=(\lambda_1,\lambda_2)$, $0\leq \lambda_1<\lambda_2\leq\infty$, such
 that
 \be
 \label{unb_lambda}
 \lambda=\int_0^{\infty} \frac{\lambda^\theta}{\Gamma(\theta)}
 x^{\theta-1} e^{-\lambda x} \delta(x) \ud x,
 \ \ \ \ \lambda\in I,
 \ee
 where
 $\delta(x):=K(\theta/n)\int_{0}^1 y^{\theta/n-1} (1-y)^{(n-3)/2}
 g(y) w(x,y) \ud y$. Certainly, the function $w$ is assumed integrable,
 that is, $\E_{(\theta,\lambda)}|w(X,Y)|<\infty$, $\lambda\in I$. Hence,
 $\delta(x)$ is a.e.\ finite, and we can redefine it, if necessary, to
 be finite for all $x>0$. Notice that the assumption
 $\E_{(\theta,\lambda)}|w(X,Y)|<\infty$
 for $\lambda\in I$ shows that
 \[
 %\label{finite}
 \int_0^\infty e^{-\lambda x} x^{\theta-1}|\delta(x)|\ud x
 \leq
 \int_0^\infty \int_0^1 e^{-\lambda x} x^{\theta-1}
 y^{\theta/n-1} (1-y)^{(n-3)/2}
 g(y) |w(x,y)| \ud y \ud x<\infty, \ \ \ \lambda\in I.
 \]
 We now rewrite (\ref{unb_lambda})
 as
 \[
 \frac{\Gamma(\theta)}{\lambda^\theta}=\frac{1}{\lambda}
 \int_0^{\infty}
  e^{-\lambda x} x^{\theta-1} \delta(x) \ud x,
 \ \ \ \lambda\in I.
 \]
 The left-hand side of the above equation equals to
 $\int_0^{\infty} e^{-\lambda x} x^{\theta-1} \ud x$, while
 its right-hand side can be written as
 \[
 \left(\int_0^{\infty}
  e^{-\lambda x_2} \ud x_2\right)
 \left(\int_0^{\infty} e^{-\lambda x_1} x_1^{\theta-1} \delta(x_1) \ud x_1\right)
 =\int_{0}^\infty  e^{-\lambda x} \left\{\int_0^{x}
 t^{\theta-1} \delta(t)\ud t\right\} \ud x,
 \]
 from Fubini's theorem. Therefore, we conclude that
 \[
 \int_0^{\infty} e^{-\lambda x} x^{\theta-1} \ud x
 =
 \int_{0}^\infty  e^{-\lambda x} \left\{\int_0^{x}
 t^{\theta-1} \delta(t)\ud t\right\} \ud x, \ \ \ \lambda\in I,
 \]
 and Lemma \ref{lem.laplace} shows that
 $x^{\theta-1}=\int_0^{x} t^{\theta-1} \delta(t)\ud t$, a.e.\ in $(0,\infty)$.
 Hence, the only possibility is $\delta(x)=(\theta-1)/x$, a.e., and
 since $\delta(x)\leq 0$, this contradicts (\ref{unb_lambda}).
 %If $\theta=1$ we must have $\delta(x)=0$, a.e.;  then,
 %$\delta(X)$ is, certainly, not unbiased.
 %If $0<\theta<1$, no $\lambda>0$ exists such that
 %the function
 %$x\mapsto  e^{-\lambda x} x^{\theta-1}|\delta(x)|=(1-\theta)
 %e^{-\lambda x}/x$ is integrable in $(0,\infty)$, and this
 %contradicts (\ref{finite}).
 Given $n=3$ observations from $\mathcal{G}(\alpha,\lambda)$,
 it can be seen (we omit the details) that the relation
 $\E_{(\alpha,\lambda)}w(X,Y)=\lambda$ for all $\alpha\in(\alpha_1,\alpha_2)$,
 $\lambda\in(\lambda_1,\lambda_2)$ [with $\alpha_1\geq 1/3$, $\lambda_1\geq 0$]
 implies $w(x,y)=x^{-1}(y|G'(y)|/G(y)-1)$, a.e.\ in $(0,\infty)^2$.
 Since $\E_{(\alpha,\lambda)}X^{-1}=\lambda/(3\alpha-1)$, this estimator
 would be unbiased for $\lambda$ if and only if
 $\E_\alpha Y|G'(Y)|/G(Y)=3\alpha$;
 however, we showed in Theorem \ref{theo.unbiased} that this is impossible.
 Hence, there is no unbiased estimator for $\lambda$ when $n=3$, and therefore,
 for $n\leq 3$, too.
 \medskip
 \end{pr}

 Several authors use the scale parameter $\beta=1/\lambda$ in defining
 the Gamma density, that is,
 \[
 f(x;\alpha,\beta)=\frac{x^{\alpha-1}}{\beta^\alpha \Gamma(\alpha)} e^{-x/\beta},
 \  \ \ x>0.
 \]
 We now proceed to verify that the function $G$ in (\ref{g_(n-1)})
 also appears
 in the form of the UMVUE of both $1/\alpha$ and $1/\lambda$.
 \begin{theo}{\rm
 \label{theo.unb.1/a}
 For $n\geq 2$, the UMVUE of $1/\alpha$ is given by
 \[
 u_2(Y):=\frac{(1-Y)^{-(n-3)/2}}{G(Y)}\int_{Y}^{1} \frac{1}{x} (1-x)^{(n-3)/2}
 G(x) \ud x,
 \]
 with $Y$ and $G$ as in Theorem \ref{theo.unbiased}; the UMVUE
 of $1/\lambda$ is given by
 \[
 u_3(X,Y):=\overline{X} u_2(Y)=X u_2(Y)/n.
 \]
 }
 \end{theo}
 \begin{pr}{Proof} It suffices to show that $u_2$ and $u_3$ are unbiased.
 For $u_2$ we have
 \[
 \E_{\alpha} u_2(Y)= K(\alpha)\int_0^1 y^{\alpha-1} \left\{
 \int_{y}^{1} \frac{1}{x} (1-x)^{(n-3)/2} G(x)
 \ud x \right\} \ud y
 = \frac{1}{\alpha} \int_0^1 f_Y(x;\alpha) \ud x =\frac{1}{\alpha},
 \]
 where the change in the order of integration is justified by Tonelli's theorem.
 The result for $u_3$ is obvious, since
 $\E_{(\alpha,\lambda)}\overline{X}=\alpha/\lambda$ and the random
 variables $\overline{X}$, $Y$ are independent.
 \end{pr}

 \begin{rem}{\rm
 \label{rem.closed.form}
 (Alternative derivation of the closed-form Ye-Chen-type
 estimators, and their asymptotic efficiency).
 For $n=2$, $u_2$ and $u_3$ admit closed forms because $G(y)\equiv 1$.
 More specifically, let $X_1,X_2$ be independently distributed $\mathcal{G}(\alpha,\lambda)$
 random variables, so that $\overline{X}=(X_1+X_2)/2$ and $Y=4 X_1 X_2/(X_1+X_2)^2$.
 Then,
 \[
 u_2(Y)=
 %\frac{(1-Y)^{-(n-3)/2}}{g(Y)}\int_{Y}^{1} \frac{1}{x} (1-x)^{(n-3)/2} g(x) \ud x
 %=
 \sqrt{1-Y}\int_{Y}^{1} \frac{1}{x\sqrt{1-x}} \ud x =
 \sqrt{1-Y} \left(2 \log\left(1 + \sqrt{1 - Y}\right) - \log Y \right).
 \]
 Therefore, after some algebra, one can verify that the unbiased estimators $u_2$
 and $u_3$ (when $n=2$) can be written as
 \begin{eqnarray*}
 &&
 u_2(Y)= T_2(X_1,X_2)=\frac{(X_2-X_1) (\log X_2-\log X_1)}{X_1+X_2},
 \\
 &&
 u_3(X,Y)=T_3(X_1,X_2)=\frac{1}{2}(X_2-X_1) (\log X_2-\log X_1).
 \end{eqnarray*}
 From the above form of $u_3=T_3$ it follows that
 $\E_{(\alpha,\lambda)} T_3(X_1,X_2)=\Cov_{(\alpha,\lambda)}(X_1,\log X_1)$.
 The well-known Stein-type identity for the Gamma distribution --
 see, e.g., Afendras {\it et al.}\ (2011) -- states
 that for any absolutely continuous function $w:(0,\infty)\to \RR$,
 \[
 \Cov_{(\alpha,\lambda)} (X_1,w(X_1)) =\frac{1}{\lambda} \E_{(\alpha,\lambda)}
 X_1 w'(X_1),
 \]
 provided $\E_{(\alpha.\lambda)}X_1|w'(X_1)|<\infty$. It is now clear that,
 if we apply Stein's identity to $w(x)=\log x$, we conclude
 the unbiasedness of $u_3$. Consider a random sample
 $X_1,\ldots,X_n$ of size $n\geq 2$ from any distribution; then,
 as is well-known, under suitable moment conditions,
 an unbiased estimator of $\Cov(X_1,w(X_1))$ is given by the
 corresponding $U$-statistic with kernel
 $K(X_1,X_2)=(X_2-X_1)(w(X_2)-w(X_1))/2$, that is,
 \[
 U_n:=\frac{2}{n(n-1)} \sum_{i<j} K(X_i,X_j)
 %=\frac{1}{n-1}\sum_{i=1}^n
 %(X_i-\overline{X})(w(X_i)-\overline{w}(X))
 =\frac{1}{n-1}\sum_{i=1}^n X_i \ w(X_i) -\frac{n}{n-1}\overline{X} \ \overline{w}(X),
 \]
 where $\overline{X}=n^{-1} \sum_{i=1}^n X_i$,
 $\overline{w}(X)=n^{-1} \sum_{i=1}^n w(X_i)$.
 With $w(x)=\log x$
 we obtain the highly efficient unbiased estimator of the scale parameter
 $1/\lambda$,
 \[
 \widehat{\lambda^{-1}}(n)=U_n=\frac{1}{n-1}\left\{
 \sum_{i=1}^n X_i \log X_i -\overline{X} \
 \sum_{i=1}^n \log X_i\right\},
 \]
 introduced by Ye and Chen (2017) (see their Theorem 4.1),
 using generalized
 likelihood equations. From the well-known asymptotic theory
 of $U$-statistics, Hoeffding (1948),
 we have $n\Var U_n\to 4 \zeta_1^2$ and
 $\sqrt{n}(U_n-\mu)\lawv\mathcal{N}(0,4\zeta_1^2)$, as $n\to\infty$, where
 $\mu=\E K(X_1,X_2)$, $\zeta_1^2=\Var \E(K(X_1,X_2)|X_2)$.
 Therefore, with $K=T_3$
 (equivalently, $w(x)=\log x$) we obtain
 \[
 \sqrt{n}\left(\widehat{\lambda^{-1}}(n)-\lambda^{-1}\right)
 \lawv\mathcal{N}\Big(0,\big(1+\alpha\psi_1(\alpha)\big)/\lambda^2\Big),
 \ \ \ \mbox {as } \ n\to\infty,
 \]
 in accordance to Theorem 3.2 in Ye and Chen (2017), since $1/\lambda=\beta$.
 Hence, the asymptotic
 efficiency of the estimator $\widehat{\lambda^{-1}}(n)$
 (compared to the MLE, see Remark \ref{rem.FI}) is given by $\operatorname{ARE}_{\widehat{\lambda^{-1}}(n)}(\alpha,\lambda)=
 \psi_1(\alpha)/(\alpha^2\psi_1(\alpha)^2-1)$, and it is independent
 of $\lambda$. Numerical calculations show that the asymptotic efficiency is greater
 than $96.4\%$ for all values of $\alpha>0$ (the minimum of the function
 $\psi_1(\alpha)/(\alpha^2\psi_1(\alpha)^2-1)$ is attained
 at $\alpha=\alpha^*= .920835\ldots$). Applying the delta-method (using
 the mapping $x\mapsto x^{-1}$) we immediately obtain
 the corresponding result for the rate parameter
 $\lambda$, namely,
 $\sqrt{n}\left(1/\widehat{\lambda^{-1}}(n)-\lambda\right)
 \lawv \mathcal{N}(0,(1+\alpha\psi_1(\alpha))\lambda^2)$, $n\to\infty$.
 Certainly, the asymptotic efficiency is the same:
 \[
 \operatorname{ARE}_{1/\widehat{\lambda^{-1}}(n)}(\alpha,\lambda)
 =\frac{\lambda^2 \psi_1(\alpha)/(\alpha\psi_1(\alpha)-1)}{(1+\alpha\psi_1(\alpha))\lambda^2}
 =\frac{\psi_1(\alpha)}{\alpha^2\psi_1(\alpha)^2-1}>0.964.
 \]

 In a similar fashion, a closed-form unbiased estimator for $1/\alpha$ can
 be constructed as an $U$-statistic with kernel $K(X_1,X_2)=T_2(X_1,X_2)$.
 However, it is difficult to evaluate the asymptotic variance of the
 resulting estimator, because $\Var \E (T_2(X_1,X_2)|X_2)$ is intractable.
 Ye and Chen proposed a different unbiased estimator for $\alpha^{-1}$,
 which is equal to the UMVUE if and only if $n=2$, namely
 %is given
 %by Theorem 4.1 in Ye and Chen (2017), namely
 \[
 \widehat{\alpha^{-1}}(n)=\frac{1}{n-1}
 \frac{\sum_{i=1}^n X_i\log X_i -\overline{X} \sum_{i=1}^n \log X_i }{\overline{X}}.
 \]
 From Theorem 3.2
 in Ye and Chen (2017) it follows that
 $\sqrt{n}\left(1/\widehat{\alpha^{-1}}(n)-\alpha\right)
 \lawv\mathcal{N}\left(0,\alpha v(\alpha)\right)$
 and
 $\sqrt{n}\left(\widehat{\alpha^{-1}}(n)-\alpha^{-1}\right)
 \lawv\mathcal{N}\left(0,\alpha^{-3}v(\alpha)\right)$,
 as $n\to\infty$,
 with $v(\alpha)=\alpha^2\psi_1(\alpha)+\alpha-1$.
 Therefore, both closed-form estimators $\widehat{\alpha^{-1}}(n)$ (for $1/\alpha$)
 and $1/\widehat{\alpha^{-1}}(n)$ (for $\alpha$) have asymptotic efficiency
 (compared to the MLE -- see Remark \ref{rem.FI}) given by
 $\rho(\alpha)
 =\big[(\alpha \psi_1(\alpha)-1)(\alpha^2\psi_1(\alpha)+\alpha-1)\big]^{-1}$.
 Numerical computations show that $\rho(\alpha)\geq \rho(\alpha^{*})=.9281\ldots$
 for all $\alpha>0$ (where $\alpha^*=.41541\ldots$), hence,
 the asymptotic efficiency is greater then $92.8\%$.
 It should be noted, however, that, due to Jensen's inequality,
 both estimators
 $1/\widehat{\lambda^{-1}}(n)$ (for $\lambda$) and
 $1/\widehat{\alpha^{-1}}(n)$ (for $\alpha$) are
 positively biased.

 The present powerful estimation procedure of Ye and Chen was further analyzed
 by Louzada {\it et al.}\  (2019),
 showing that the bias-corrected estimator
 \[
 \widehat{\alpha}_n=
 \frac{n-3}{(n-1)\widehat{\alpha^{-1}}(n)}
 +\frac{2}{5\left(n+(n-1)\widehat{\alpha^{-1}}(n)\right)}
 \
 \Bigg(
 \frac{11}{3}-
 \frac{2n}{n+(n-1)\widehat{\alpha^{-1}}(n)}
 \Bigg)
 \]
 improves upon the estimator $1/\widehat{\alpha^{-1}}(n)$ of $\alpha$
 for small values of $n$.
 }
 \end{rem}

 \section{Closed form estimators for the two-parameter Beta distribution}
 \label{sec.3}
 Through the present section we consider a random sample $X_1,\ldots,X_n$
 ($n\geq 2$) from the
 %$B(\alpha,\beta)$
 Beta
 distribution with density as in (\ref{beta.distribution}),
 %\[
 %f(x)=\frac{1}{B(\alpha,\beta)} x^{\alpha-1}(1-x)^{\beta-1}, \ \ \ 0<x<1,
 %\]
 where the parameter ${\bbb \theta}=(\alpha,\beta)\in(0,\infty)^2$.
 %, and
 %$B(\alpha,\beta)=\Gamma(\alpha)\Gamma(\beta)/\Gamma(\alpha+\beta)$.
 The Fisher information matrix (based on a single observation)
 is given by
 \begin{eqnarray*}
 &&
 J=J(\alpha,\beta)=
 \left(
 \begin{array}{cc}
 \psi_1(\alpha)-\psi_1(\alpha+\beta) & -\psi_1(\alpha+\beta) \\
 -\psi_1(\alpha+\beta) & \psi_1(\beta)-\psi_1(\alpha+\beta)
 \end{array}
 \right),
 \mbox{ with inverse }
 \\
 &&
 J^{-1}=\frac{1}{\psi_1(\alpha)\psi_1(\beta)-\psi_1(\alpha+\beta)
 (\psi_1(\alpha)+\psi_1(\beta))}
 \left(
 \begin{array}{cc}
 \psi_1(\beta)-\psi_1(\alpha+\beta) & \psi_1(\alpha+\beta) \\
 \psi_1(\alpha+\beta) & \psi_1(\alpha)-\psi_1(\alpha+\beta)
 \end{array}
 \right),
 \end{eqnarray*}
 with $\psi_1(x)=\left(\log \Gamma(x)\right)''$, as in Section \ref{sec.2}.
 The complete sufficient statistic has now the form $(X,Y)
 =\left(\prod_{i=1}^n X_i,\prod_{i=1}^n (1-X_i)\right)$, and an extra difficulty,
 compared to the Gamma-case,
 is due to the lack of independence of $X$ and $Y$.

 The proposed Ye-Chen-type (biased) estimators for $\alpha$ and $\beta$
 are as follows:
 \begin{eqnarray}
 \label{unb.a}
 \widehat{a}_n=\frac{n\overline{X}_n}{\sum_{i=1}^n X_i \log(X_i/(1-X_i))-
 \overline{X}_n \sum_{i=1}^n \log(X_i/(1-X_i))},
 \\
 \label{unb.b}
 \widehat{\beta}_n=\frac{n(1-\overline{X}_n)}{\sum_{i=1}^n X_i \log(X_i/(1-X_i))-
 \overline{X}_n \sum_{i=1}^n \log(X_i/(1-X_i))},
 \end{eqnarray}
 where $\overline{X}_n=n^{-1}\sum_{i=1}^n X_i$.
 In order to make clear how these estimators have been derived, we recall
 the Stein-type covariance identity for the Beta distribution,
 \[
  \Cov(X,w(X))=\frac{1}{\alpha+\beta}\E X(1-X) w'(X),
 \]
 which holds true for sufficiently smooth functions, namely, when $w:(0,1)\to\RR$
 is absolutely continuous and its derivative, $w'$, satisfies $\E X(1-X)|w'(X)|<\infty$;
 see, e.g., Afendras {\it et al.}\ (2011). The special choice
 $w(x)=\log(x/(1-x))$ leads to $\Cov(X,w(X))=1/(\alpha+\beta)$. Therefore,
 the $U$-statistic $U_n$ with kernel
 \[
 K(X_1,X_2)=\frac{1}{2}(X_2-X_1)\log\frac{X_2(1-X_1)}{X_1(1-X_2)}
 \]
 is an unbiased estimator of the parameter $1/(\alpha+\beta)$. One finds
 \[
  U_n =\frac{2\sum_{i<j} K(X_i,X_j)}{n(n-1)} =
  \frac{1}{n-1}\left(
  \sum_{i=1}^n X_i \log(X_i/(1-X_i))-
 \overline{X}_n \sum_{i=1}^n \log(X_i/(1-X_i))
 \right)=\frac{n V_n}{n-1},
 \]
 say. From the SLLN,
 $V_n\to \Cov(X,\log(X/(1-X)))=1/(\alpha+\beta)$, a.s., as $n\to\infty$; thus,
 $1/V_n\to \alpha+\beta$, a.s. Since, obviously,
 $\overline{X}_n\to \alpha/(\alpha+\beta)$ and $1-\overline{X}_n\to \beta/(\alpha+\beta)$,
 we conclude that $\widehat{\alpha}_n=\overline{X}_n/V_n\to \alpha$
 and $\widehat{\beta}_n=(1-\overline{X}_n)/V_n\to \beta$, a.s., so both estimators
 are strongly consistent.
 In order to investigate their asymptotic efficiency, we have to calculate
 explicitly the variance of their asymptotic normal distribution.
 This is the subject of the following
 \begin{theo}{\rm
 \label{theo.as.normal}
 The joint limiting distribution of $\widehat{\alpha}_n$, $\widehat{\beta}_n$
 (see (\ref{unb.a}), (\ref{unb.b})) is normal. More precisely, as $n\to\infty$,
 \be
 \label{normal.limit.a.b}
 \sqrt{n}\left(
 \begin{array}{c}
 \widehat{\alpha}_n-\alpha \\
 \widehat{\beta}_n-\beta
 \end{array}\right)
 \lawv\mathcal{N}({\bbb 0}, \Delta), \ \ \
 \mbox{where} \ \ \
 \Delta=\frac{1}{\alpha+\beta+1}\left(
 \begin{array}{cc}
 v_{11}(\alpha,\beta) & v_{12}(\alpha,\beta) \\
 v_{12}(\alpha,\beta) & v_{22}(\alpha,\beta)
 \end{array}
 \right)
 \ee
 with
 \begin{eqnarray*}
   v_{11}(\alpha,\beta)
  \hspace*{-1ex}&=&\hspace*{-1ex}
   \alpha \left( \alpha (\alpha + \beta+1) + \alpha^2 \beta
   \
   \Big[\psi_1(\alpha) + \psi_1(\beta)\Big]-\beta \right),
   \\
   v_{12}(\alpha,\beta)
   \hspace*{-1ex}&=&\hspace*{-1ex}
   \alpha^2\beta+\alpha\beta^2 +  \alpha^2 \beta^2 \
   \Big[\psi_1(\alpha) + \psi_1(\beta)\Big]
   -\alpha^2-\beta^2,
   \\
   v_{22}(\alpha,\beta)
   \hspace*{-1ex}&=&\hspace*{-1ex}
   \beta \left(\beta (\alpha + \beta+1) + \alpha \beta^2
   \
   \Big[\psi_1(\alpha) + \psi_1(\beta)\Big]-\alpha\right).
  \end{eqnarray*}
  }
  \end{theo}
  \begin{pr}{Proof} The result is quite simple,
  since it is a by-product of the CLT and the delta-method.
  However, the calculations are somewhat involved.
  We define the random variables $Y_i=\log(X_i/(1-X_i))$
  and $Z_i=X_i\log(X_i/(1-X_i))$. It follows after some
  algebra that the mean vector and the dispersion matrix
  of $(X_i,Y_i,Z_i)$, $i=1,\ldots,n$, are given by
  \[
  {\bbb \mu}=
  \left(
  \begin{array}{c}
  \mu_1
  \\
  \mu_2
  \\
  \mu_3
  \end{array}
   \right)
   =
  \frac{1}{\alpha+\beta}
  \left(
  \begin{array}{c}
  \alpha
  \\
  (\alpha+\beta)(\psi(\alpha)-\psi(\beta))
  \\
  1+\alpha(\psi(\alpha)-\psi(\beta))
  \end{array}
   \right),
   \ \ \ \ \ \
  \Sigma=\left(
  \begin{array}{ccc}
  \sigma_{11} &\sigma_{12}& \sigma_{13}
  \\
  \sigma_{12} &\sigma_{22}& \sigma_{23}
  \\
  \sigma_{13} &\sigma_{23}& \sigma_{33}
  \end{array}
   \right),
  \]
  where, with $\Psi=\psi(\alpha)-\psi(\beta)$ (recall $\psi(x)=\left(\log \Gamma(x)\right)'$)
  and $\Psi_1=\psi_1(\alpha)+\psi_1(\beta)$,
  \begin{eqnarray*}
   \sigma_{11}
   \hspace*{-1ex}&=&\hspace*{-1ex}
   \frac{\alpha\beta}{(\alpha+\beta)^2(\alpha + \beta+1)},
    \\
   \sigma_{12}
   \hspace*{-1ex}&=&\hspace*{-1ex}
   \frac{1}{\alpha+\beta},
   \\
   \sigma_{13}
   \hspace*{-1ex}&=&\hspace*{-1ex}
   \frac{\beta + \alpha (\alpha + \beta)+\alpha\beta \ \Psi}{(\alpha+\beta)^2(\alpha+\beta+1)},
   \\
   \sigma_{22}
  \hspace*{-1ex}&=&\hspace*{-1ex}
   \Psi_1, \begin{array}{c} \\ \frac{}{} \end{array}
   \\
   \sigma_{23}
  \hspace*{-1ex}&=&\hspace*{-1ex}
   \frac{\Psi+\alpha \Psi_1}{\alpha+\beta},
   \\
   \sigma_{33}
  \hspace*{-1ex}&=&\hspace*{-1ex}
  \frac{\alpha\beta \ \Psi^2+2(\alpha^2+\alpha\beta+\beta)
  \Psi+\alpha(\alpha+1)(\alpha+\beta)\Psi_1
  +\alpha+\beta-1}{(\alpha+\beta)^2(\alpha + \beta+1)}.
  \end{eqnarray*}
  Hence, from the classical CLT,
  \[
  \sqrt{n}\left(
 \begin{array}{c}
 \overline{X}_n-\mu_1
 \\
 \overline{Y}_n-\mu_2
 \\
 \overline{Z}_n-\mu_3
 \end{array}\right)
  \lawv\mathcal{N}({\bbb 0}, \Sigma), \ \ \ \mbox{as } \ n\to\infty,
  \]
  where $\overline{X}_n=n^{-1}\sum_{i=1}^n X_i$, \ $\overline{Y}_n=n^{-1}\sum_{i=1}^n Y_i$,
  \
  $\overline{Z}_n=n^{-1}\sum_{i=1}^n Z_i$.
  Consider the functions $k_1(x,y,z)=x/(z-xy)$, $k_2(x,y,z)=(1-x)/(z-xy)$,
  so that $k_1(\overline{X}_n,\overline{Y}_n,\overline{Z}_n)=\widehat{\alpha}_n$ and
  $k_2(\overline{X}_n,\overline{Y}_n,\overline{Z}_n)=\widehat{\beta}_n$.
  The Jacobian matrix is
  \[
  A(x,y,z):=
  \left(
 \begin{array}{ccc}
  \frac{\partial}{\partial x} k_1(x,y,z) & \frac{\partial}{\partial y} k_1(x,y,z)
  & \frac{\partial}{\partial z} k_1(x,y,z)
  \\
  \frac{\partial}{\partial x} k_2(x,y,z) & \frac{\partial}{\partial y} k_2(x,y,z)
  & \frac{\partial}{\partial z} k_2(x,y,z)
  \end{array}\right)
  =\frac{1}{(xy-z)^2}
  \left(
 \begin{array}{ccc}
  z & x^2  & -x
  \\
  y-z & x-x^2 & x-1
  \end{array}\right).
  \]
  Hence, with $A_0:=A(\mu_1,\mu_2,\mu_3)$, a straightforward application
  of the delta-method yields
  \[
  \sqrt{n}\left(
  \begin{array}{c}
  k_1(\overline{X}_n,\overline{Y}_n,\overline{Z}_n)-k_1(\mu_1,\mu_2,\mu_3)
  \\
  k_2(\overline{X}_n,\overline{Y}_n,\overline{Z}_n)-k_2(\mu_1,\mu_2,\mu_3)
  \end{array}\right)
  =
  \sqrt{n}\left(
  \begin{array}{c}
  \widehat{\alpha}_n-\alpha
  \\
  \widehat{\beta}_n-\beta
  \end{array}\right)
   \lawv\mathcal{N}({\bbb 0}, A_0\Sigma A_0^T),
   \ \ \ \mbox{as } \ \ n\to\infty,
  \]
  because, obviously,
  $k_1(\mu_1,\mu_2,\mu_3)=\alpha$ and $k_2(\mu_1,\mu_2,\mu_3)=\beta$.
  Finally, a lengthy computation reveals that
  $A_0\Sigma A_0^T=\Delta$ as in (\ref{normal.limit.a.b}),
  completing the proof.
  \medskip
  \end{pr}

  \noindent
  For any unbiased estimator $U$ of $h(\alpha,\beta)$ (provided that
  such an estimator exists),
  the CR-bound yields the inequality
  $n \Var U \geq \nabla h(\alpha,\beta)^T J^{-1}
  \nabla h(\alpha,\beta)$; in particular, with $h(\alpha,\beta)=\alpha$
  we get
  \[
  n\Var U \geq  \frac{\psi_1(\beta)-\psi_1(\alpha+\beta)}{\psi_1(\alpha)\psi_1(\beta)-\psi_1(\alpha+\beta)
  (\psi_1(\alpha)+\psi_1(\beta))},
  \]
  and the right-hand side is the asymptotic variance of the MLE of $\alpha$,
  which, as is well known, does not admit a closed form.
  For $h(\alpha,\beta)=\beta$ we obtain a similar inequality,
  interchanging the roles of $\alpha$ and $\beta$. Hence, from
  Theorem \ref{theo.as.normal},
  the asymptotic efficiency $\rho_1=\rho_1(\alpha,\beta)$
  of $\widehat{\alpha}_n$ is given by
  \[
  \rho_1=
  \frac{(\alpha+\beta+1)(\psi_1(\beta)-\psi_1(\alpha+\beta))/\alpha}
  {\Big( \alpha (\alpha + \beta+1) + \alpha^2 \beta
   \
   \big[\psi_1(\alpha) + \psi_1(\beta)\big]-\beta \Big)
   \Big(\psi_1(\alpha)\psi_1(\beta)-\psi_1(\alpha+\beta)
   (\psi_1(\alpha)+\psi_1(\beta))\Big)},
  \]
  while the asymptotic efficiency of $\widehat{\beta}_n$ equals to
  $\rho_2(\alpha,\beta)=\rho_1(\beta,\alpha)$.
  Despite the fact that the function $\rho_1$ is too complicated to enable
  a complete study, we can report some properties supported by
  numerical evidence. For fixed $\beta>0$, define $\phi(\alpha)=\rho_1(\alpha,\beta)$.
  The function $\phi$ has limits $1$, as $\alpha\to 0+$, and
  $\psi_1(\beta)/(\beta^2\psi_1(\beta)^2-1)>.964$, as $\alpha\to+\infty$;
  see Remark \ref{rem.closed.form}. Hence, either the infimum
  equals to the limit at $+\infty$, or it is attained at some point
  $\alpha^*=\alpha^*(\beta)$ such that $\phi'(\alpha^*)=0$.
  Several plots showed that the graph of $\phi$ (for any fixed
  $\beta$) first decreases, then increases, and finally decreases.
  We calculated numerically the smallest  stationary point, $\alpha^*$, of $\phi$
  (always lying between $0$ and $1$), resulting
  to the  inequality $\phi(\alpha)\geq\min\{\phi(\alpha^*),.964\}$.
  Numerical values of $\rho_1^*(\beta)=\rho_1(\alpha^*(\beta),\beta)$, as well as
  minimization points $\alpha^*(\beta)$,
  are presented in Table 1. These values indicate
  that $\rho_1(\alpha^*(\beta),\beta)\to\min_{x>0}
  \big[(x\psi_1(x)-1)(x^2\psi_1(x)+x-1)\big]^{-1}=.9281\ldots$,
  and $\alpha^*(\beta)\to \alpha^*=.41541\ldots$ ($\alpha^*$ is the same
  as in Remark
  \ref{rem.closed.form}), as $\beta\to\infty$;
  however, we was not able to provide a proof.

  The values of Table 1
  indicate that for almost the entire parametric range,
  the asymptotic efficiency of $\widehat{\alpha}_n$ is no less that $90\%$.
  The same is clearly true for $\widehat{\beta}_n$, since
  $\rho_2(\alpha,\beta)=\rho_1(\beta,\alpha)$. This verifies
  that the proposed estimating procedure is quite powerful and, at the
  same time, it produces simple, closed-form, estimators.

  \newpage
  \noindent
  \centerline{\bf Table 1.}

  \noindent
  \centerline{Least asymptotic efficiency
  $\rho_1^*(\beta)=\min_{\alpha>0}\rho_1(\alpha,\beta)=\rho_1(\alpha^*(\beta),\beta)$,
  for $\beta=.05(.05)10$.}
  \medskip

  \noindent
  {\footnotesize
  \begin{tabular}{|lll|lll|lll|lll|}
  \hline
  $\beta$ & $\alpha^*(\beta)$ & $\rho_1^*(\beta)$ &
  $\beta$ & $\alpha^*(\beta)$ & $\rho_1^*(\beta)$ &
  $\beta$ & $\alpha^*(\beta)$ & $\rho_1^*(\beta)$ &
  $\beta$ & $\alpha^*(\beta)$ & $\rho_1^*(\beta)$ \\
  \hline
     .05 &  .179398 &  .712626  &  .1 &  .229471 &  .784172  &  .15 &   .258632 &   .828715 &
   .2 &   .277437 &   .85901  \\
   .25 &   .290137 &    .880541 &  .3 &   .298956 &   .89629 &  .35 &   .305205 &   .908065 &
   .4 &   .309705 &   .917024 \\
   .45 &   .312992 &   .923944 &  .5 &   .315422 &   .929356 &  .55 &   .317239 &    .933636 &
   .6 &   .318612 &   .937052 \\
   .65 &   .319661 &    .939802 & .7 &   .320472 &   .942029 &  .75 &   .321107 &    .943844 &
   .8 &   .321614 &   .945328 \\
   .85 &   .322026 &    .946546 &  .9 &   .322371 &   .947547 &  .95 &   .322669 &    .94837 &
  1. &   .322936 &   .949046   \\
  1.05 &   .323182 &   .9496 &  1.1 &    .323418 &   .950051 &  1.15 &   .32365 &   .950416 &
  1.2 &   .323884 &    .950709 \\
  1.25 &   .324122 &   .95094 &  1.3 &   .324369 &    .951118 &   1.35 &   .324626 &   .951251 &
  1.4 &   .324894 &    .951345 \\
  1.45 &   .325174 &   .951405 &  1.5 &   .325466 &  .951437 &  1.55 &   .325771 &   .951445 &    1.6 &   .326088 &   .95143  \\
  1.65 &   .326417 &   .951397 &  1.7 &   .326757 &  .951348 &    1.75 &   .327109 &   .951284 &    1.8 &   .327471 &   .951209 \\
  1.85 &   .327842 &   .951122 &    1.9 &   .328222 &   .951027 &    1.95 &   .32861 &  .950924 &   2. &   .329006 &   .950814 \\
  2.05 &   .329409 &   .950697 &   2.1 &   .329818 &  .950576 &  2.15 &   .330232 &   .950451 &    2.2 &   .330652 &   .950321 \\
  2.25 &   .331075 &   .950189 &   2.3 &   .331502 &  .950054 &    2.35 &   .331933 &   .949917 &    2.4 &   .332366 &   .949778 \\
  2.45 &   .332801 &   .949638 &    2.5 &   .333239 &   .949496 &    2.55 &   .333677 & .949354 &   2.6 &   .334117 &  .949212  \\
  2.65 &   .334557 &   .949069 &    2.7 &   .334997 & .948926 &  2.75 &   .335438 &   .948783 &    2.8 &   .335878 &  .94864 \\
  2.85 &   .336317 &   .948497 &  2.9 &   .336756 &  .948355 &    2.95 &   .337194 &   .948214 &    3. &   .33763 &   .948073 \\
  3.05 &   .338065 &   .947933 &  3.1 &   .338498 &   .947794 &    3.15 &   .33893 &   .947656 &   3.2 &   .339359 &   .947519 \\
  3.25 &   .339786 &   .947383 &  3.3 &   .340212 &  .947248 &  3.35 &   .340635 &   .947114 &    3.4 &   .341055 &   .946981 \\
  3.45 &   .341473 &   .946849 &  3.5 &   .341888 &  .946718 &    3.55 &   .342301 &   .946589 &    3.6 &   .342711 &   .946461 \\
  3.65 &   .343118 &   .946334 &  3.7 &   .343522 &  .946208 &   3.75 &   .343924 &   .946084 & 3.8 &   .344322 &   .945961 \\
  3.85 &   .344718 &   .945839 &  3.9 &   .34511 &   .945718 &  3.95 &   .3455 &   .945599 &
  4. &   .345886 &   .945481 \\
  4.05 &   .34627 &   .945364 &  4.1 &   .34665 &  .945248 &   4.15 &   .347027 &   .945134 &    4.2 &   .347401 &   .945021 \\
  4.25 &   .347772 &   .944909 &  4.3 &   .34814 &  .944798 &    4.35 &   .348505 &   .944689 & 4.4 &   .348867 &   .94458 \\
  4.45 &   .349226 &   .944473 &  4.5 &   .349581 &   .944367 &   4.55 &   .349934 &   .944262 &    4.6 &   .350284 &   .944159 \\
  4.65 &   .35063 &   .944056 &  4.7 &   .350974 &   .943955 &    4.75 &   .351314 &   .943855 &    4.8 &   .351652 &  .943756 \\
  4.85 &   .351986 &   .943658 &  4.9 &   .352318 &   .943561 &  4.95 &  .352647 &   .943465 &
  5. &   .352973 &  .94337 \\
  5.05 &   .353295 &   .943276 &  5.1 &   .353616 &  .943184 &  5.15 &   .353933 &   .943092 &   5.2 &   .354247 &  .943001 \\
  5.25 &   .354559 &   .942911 &  5.3 &   .354868 &   .942823 &    5.35 &   .355174 &   .942735 &    5.4 &   .355477 &  .942648 \\
  5.45 &   .355778 &   .942562 &  5.5 &   .356076 &  .942477 &    5.55 &   .356371 &   .942393 &   5.6 &   .356664 &  .94231 \\
  5.65 &   .356955 &   .942228 &  5.7 &   .357242 &   .942146 &   5.75 &   .357527 &   .942066 &    5.8 &   .35781 &  .941986 \\
  5.85 &   .35809 &   .941907 &  5.9 &   .358368 &  .941829 &   5.95 &   .358643 &   .941752 &    6. &  .358916 &  .941676 \\
  6.05 &   .359186 &   .9416 &  6.1 &   .359454 &  .941525 &    6.15 &   .35972 &   .941451 &   6.2 &   .359983 &  .941378 \\
  6.25 &   .360244 &   .941305 &  6.3 &   .360503 &  .941234 &  6.35 &   .36076 &   .941163 &    6.4 &   .361014 &   .941092 \\
  6.45 &   .361266 &   .941023 &  6.5 &   .361516 &  .940954 &   6.55 &   .361764 &   .940886 &    6.6 &   .36201 &   .940818 \\
  6.65 &   .362254 &   .940751 &  6.7 &   .362495 &   .940685 &    6.75 &   .362735 &   .94062 &  6.8 &   .362972 &   .940555 \\
  6.85 &   .363208 &   .940491 &  6.9 &   .363441 &   .940427 &  6.95 &   .363673 &   .940364 &    7. &   .363903 &   .940302 \\
  7.05 &   .36413 &   .94024 &  7.1 &   .364356 &   .940179 &    7.15 &   .36458 &   .940118 &    7.2 &   .364802 &   .940058 \\
  7.25 &   .365022 &   .939999 &  7.3 &   .365241 &   .93994 &   7.35 &   .365457 &   .939882 &   7.4 &   .365672 &   .939824 \\
  7.45 &   .365885 &   .939767 &  7.5 &   .366097 &  .93971 &  7.55 &   .366307 &   .939654 &    7.6 &   .366514 &   .939599 \\
  7.65 &   .366721 &   .939544 &  7.7 &   .366925 &   .939489 &    7.75 &   .367128 &   .939435 &    7.8 &   .36733 &   .939381 \\
  7.85 &   .36753 &   .939328 &  7.9 &   .367728 &   .939276 &    7.95 &   .367924 &   .939224 &  8. &   .368119 &  .939172 \\
  8.05 &   .368313 &   .939121 &  8.1 &   .368505 &   .93907 &  8.15 &   .368695 &   .93902 &    8.2 &   .368884 &   .93897 \\
  8.25 &   .369072 &   .93892 &  8.3 &   .369258 &  .938872 &   8.35 &   .369443 &   .938823 &    8.4 &   .369626 &   .938775 \\
  8.45 &   .369808 &   .938727 &  8.5 &   .369988 &   .93868 &    8.55 &   .370167 &   .938633 &   8.6 &   .370345 &   .938587 \\
  8.65 &   .370521 &   .938541 &  8.7 &   .370696 &  .938495 &   8.75 &   .37087 &   .93845 &    8.8 &   .371042 &   .938405 \\
  8.85 &   .371214 &   .93836 &  8.9 &   .371383 &  .938316 &   8.95 &   .371552 &   .938272 &    9. &   .371719 &   .938229 \\
  9.05 &   .371885 &   .938186 &  9.1 &   .37205 &  .938143 &    9.15 &   .372214 &   .938101&   9.2 &   .372376 &   .938059 \\
  9.25 &   .372538 &   .938017 &  9.3 &   .372698 &   .937976 &   9.35 &   .372857 &   .937935 &   9.4 &   .373014 &   .937894 \\
  9.45 &   .373171 &   .937854 &  9.5 &   .373327 &   .937814 &   9.55 &   .373481 &   .937774 &    9.6 &   .373634 &   .937735 \\
  9.65 &   .373786 &   .937696 &  9.7 &   .373938 &   .937657 &    9.75 &   .374088 &   .937619 &  9.8 &   .374237 &   .937581 \\
  9.85 &   .374385 &   .937543 &  9.9 &   .374531 &   .937505 &   9.95 &   .374677 &   .937468 &  10. &   .374822 &     .937431 \\
  \hline
  \end{tabular}
  }

  \begin{rem}{\rm
  \label{rem.positivity}
  (a)
  The estimators $\widehat{\alpha}_n$ and $\widehat{\beta}_n$ are positive (w.p.\ 1)
  for all $n\geq 2$. This is so because the denominator is a positive multiple
  of an $U$-statistic with symmetric kernel $K(X_1,X_2)>0$, w.p.\ $1$.

  \noindent
  (b)
  The asymptotic covariance of $\widehat{\alpha}_n$ and $\widehat{\beta}_n$
  is positive. This follows from the inequality $\psi_1(x)>x^{-2}>x^{-2}-x^{-1}$, $x>0$,
  which implies that $p(x):=\psi_1(x)-(x^{-2}-x^{-1})>0$.
  Then, from (\ref{normal.limit.a.b}),
  $v_{12}(\alpha,\beta)
  %=
   %\alpha^2\beta+\alpha\beta^2 +  \alpha^2 \beta^2 \
   %\Big[\psi_1(\alpha) + \psi_1(\beta)\Big]
   %-\alpha^2-\beta^2
   =\alpha^2\beta^2 \ \big[p(\alpha) +p(\beta) \big]>0$.
  This fact is somewhat counterintuitive, since $\widehat{\alpha}_n=\overline{X}_n/V_n$,
  $\widehat{\beta}_n=(1-\overline{X}_n)/V_n$, $V_n>0$, and the random variables
  $\overline{X}_n$ and $1-\overline{X}_n$ are, obviously, negatively correlated.

  \noindent
  (c) The limiting distribution of the estimator $U_n$
  (which is unbiased for $1/(\alpha+\beta)$) can be derived from
  Theorem \ref{theo.as.normal}, since $U_n=(n/(n-1)) V_n$ and thus,
  $\sqrt{n}(U_n-V_n)=o_p(1)$ and $\sqrt{n}(1/U_n-1/V_n)=o_p(1)$, $n\to\infty$.
  Applying the delta-method to (\ref{normal.limit.a.b}) with $k(x,y)=1/(x+y)$
  we obtain
  $\sqrt{n}\left(U_n-(\alpha+\beta)^{-1}\right)
  =\sqrt{n}\left(V_n-(\alpha+\beta)^{-1}\right)+o_p(1)
  \lawv\mathcal{N}\left(0,v^2\right)$, as $n\to\infty$,
   where
  \[
  v^2=\frac{v_{11}(\alpha,\beta)
  %k_1(\alpha,\beta)^2
  +v_{22}(\alpha,\beta)
  %k_2(\alpha,\beta)^2
  +2 v_{12}(\alpha,\beta)
  %k_1(\alpha,\beta) k_2(\alpha,\beta)
  }
  {(\alpha+\beta)^4(\alpha+\beta+1)}
  =\frac{\alpha+\beta+\alpha\beta
  \
  \big[\psi_1(\alpha) + \psi_1(\beta)\big]-1}
  {(\alpha+\beta)^2(\alpha+\beta+1)}.
  \]
  Consequently,
  $
  \sqrt{n} \ \Big(1/V_n-(\alpha+\beta)\Big)
  \lawv\mathcal{N}\left(0,(\alpha+\beta)^4 v^2\right)$, as $ n\to\infty$.
  Using these formulas one can verify that the estimators $U_n$ and $1/V_n$
  exhibit high asymptotic efficiency.

  \noindent
  (d) For small sample sizes, bias-corrected estimators
  in the sense of Louzada {\it et al.}\
  may reduce the mean squared error; this, however,
  is beyond the scope of the present work.
  %with $k_1(x,y)=k_2(x,y)=-(x+y)^{-2}$.
  }
  \end{rem}

 \end{document}